\documentclass[reqno,10pt]{amsart}

\usepackage{xcolor}
\usepackage{hyperref}
\usepackage{amssymb}
\usepackage{comment}
\usepackage{bbm}
\usepackage[all,cmtip]{xy}
\newcounter{mt}
\newcounter{mq}

\newtheorem{MainTheorem}[mt]{Theorem}
\newtheorem{Proposition}{Proposition}[section]
\newtheorem{Definition}[Proposition]{Definition}
\newtheorem{Lemma}[Proposition]{Lemma}
\newtheorem{Theorem}[Proposition]{Theorem}
\newtheorem{Corollary}[Proposition]{Corollary}
\newtheorem{Remark}[Proposition]{Remark}
\newtheorem{Conjecture}[Proposition]{Conjecture}
\newtheorem{Question}[mq]{Question}
\DeclareMathOperator{\Val}{Val}

\DeclareMathOperator{\Gr}{Gr}

\DeclareMathOperator{\Jac}{Jac}

\DeclareMathOperator{\vol}{vol}

\DeclareMathOperator{\Dens}{Dens}
\DeclareMathOperator{\Span}{Span}

\DeclareMathOperator{\Image}{Im}
\DeclareMathOperator{\Supp}{Supp}

\DeclareMathOperator{\End}{End}

\DeclareMathOperator{\Ker}{Ker}
\DeclareMathOperator{\Kl}{Kl}

\DeclareMathOperator{\Sc}{Sc}

\DeclareMathOperator{\GL}{GL}

\DeclareMathOperator{\OO}{O}
\DeclareMathOperator{\SO}{SO}

\DeclareMathOperator{\Cos}{\mathcal C}
\newcommand{\R}{\mathbb{R}}
\newcommand{\C}{\mathbb{C}}

\makeatletter
\def\moverlay{\mathpalette\mov@rlay}
\def\mov@rlay#1#2{\leavevmode\vtop{%
		\baselineskip\z@skip \lineskiplimit-\maxdimen
		\ialign{\hfil$\m@th#1##$\hfil\cr#2\crcr}}}
\newcommand{\charfusion}[3][\mathord]{
	#1{\ifx#1\mathop\vphantom{#2}\fi
		\mathpalette\mov@rlay{#2\cr#3}
	}
	\ifx#1\mathop\expandafter\displaylimits\fi}
\makeatother

\title[Quasianalyticity, uncertainty, and integral transforms]{Quasianalyticity, uncertainty, and integral transforms on higher grassmannians}
\author{Dmitry Faifman}
	\email{faifmand@tauex.tau.ac.il} 
\address{School of Mathematical Sciences, Tel Aviv University, Tel Aviv 6997801, Israel}

\usepackage[abbrev]{amsrefs}

\BibSpec{article}{%
	+{}{\PrintAuthors}  		{author}
	+{,}{ \textit}     		{title}
	+{,}{ }             		{journal}
	+{}{ \textbf}       		{volume}
	+{}{ \parenthesize} 		{date}
	+{}{, no. } 		{number}
	+{,}{ }      	      		{conference}
	+{,}{ }      	      		{book}
	+{,}{ }            		{pages}
	+{,}{ }            	 	{note}
	+{,}{ }            	 	{status}
	+{,}{  \texttt } {eprint}
	+{.}{}              {transition}
}

\begin{document}
	
	\begin{abstract}
We investigate the support of a distribution $f$ on the real grassmannian $\mathrm{Gr}_k(\mathbb R^n)$ whose spectrum, namely its nontrivial $\mathrm O(n)$-components, is restricted to a subset $\Lambda$ of all $\mathrm O(n)$-types. We prove that unless $\Lambda$ is co-sparse, $f$ cannot be supported at a point.
We utilize this uncertainty principle to prove that if $2\leq k\leq n-2$, then the cosine transform of a distribution on the grassmannian cannot be supported inside any single open Schubert cell $\Sigma^k$. The same holds for certain more general $\alpha$-cosine transforms and for the Radon transform between grassmannians, and more generally for various $\GL_n(\R)$-modules. These results are then applied to convex geometry and geometric tomography, where sharper versions of the Aleksandrov projection theorem, Funk section theorem, and Klain's and Schneider's injectivity theorems for convex valuations are obtained.

	\end{abstract}

\thanks{  
	\\\indent This research was supported by the ISRAEL SCIENCE FOUNDATION grant No. 1750/20.
}

\maketitle

\section{Overview}

\subsection{Introduction and motivation}
A quasianalytic class of functions is any class of functions on a fixed domain such that, whenever two functions coincide locally, they must coincide globally. The most common formal notion of quasianalyticity takes local coincidence to mean the coincidence of jets of the two functions at a point. This notion traces back to a question of Hadamard \cite{hadamard_question}, answered by Denjoy  and Carleman. A weaker notion of quasianalyticity goes back to S. Bernstein \cite{bernstein_quasianalytic}. A class of functions is \emph{quasianalytic in the sense of Bernstein} if, whenever two functions in this class coincide on a set with non-empty interior, they coincide globally. The latter notion has the further advantage that it allows families of non-smooth functions, or even distributions. 

Kazdan has posed in \cite{kazdan} the \emph{general quasianalytic problem} as follows. Given a class $A(\Omega)$ of smooth function on a domain $\Omega$, find subdomains $\Omega_1\subset\Omega$ such that if $f\in A(\Omega)$ vanishes on $\Omega_1$ with some derivatives, then $f$ vanishes on $\Omega$. The question admits natural variations corresponding to the various definitions of quasianalyticity. For instance, Kazdan's general quasianalytic problem in Bernstein's sense would read as follows:  Given a class $A(\Omega)$ of functions on a domain $\Omega$, find subdomains $\Omega_1\subset\Omega$ such that if $f\in A(\Omega)$ vanishes on a neighborhood of $\Omega_1$, then $f$ vanishes on $\Omega$. 

The uniqueness theorem for analytic functions is of course the first and main motivating example. Generalizing in one direction, classes of quasianalytic functions imposing restrictions on the growth rate of the Taylor series coefficients were introduced by Denjoy and Carleman \cite{denjoy_quasianalytic, carleman_quasianalytic}, and subsequently extensively studied. Bernstein in \cite{bernstein_quasianalytic} introduced a quasianalytic class by imposing restrictions on approximability by polynomials. Another direction of fundamental importance is the uniqueness problem for differential operators, with Holmgren's theorem \cite{holmgren} a prototypical example. 

This work belongs in a neighborhood of the uniqueness problem for differential operators, in that some of the classes of functions we consider are given as the image of integral transforms, which in some cases admit a PDE description, or  belong in the kernel of a differential operator. 
However, our methods apply in greater generality to certain classes of functions appearing as spaces of global sections of $\GL_n(\R)$-equivariant line bundles over the real grassmannian. As it happens, the quasianalytic property is closely linked to a new uncertainty principle on the grassmannian, which we now proceed to discuss.

\medskip 

Denote by $\Gr_k(\R^n)$ the grassmannian of $k$-dimensional linear subspaces in $\R^n$. Recall that $L^2(\Gr_k(\R^n))$ decomposes under the action of $\OO(n)$ into the direct sum of finite dimensional irreducible representations of $\OO(n)$ with multiplicity one, which we call $\OO(n)$-types. The set of all types is denoted $\Lambda_{k}^n$. For a function or distribution $f$ on $\Gr_k(\R^n)$, let $\widehat f(\lambda)$ denote its component in $\lambda\in \Lambda_k^n$. The spectrum of $f$ consists of the non-trivial components. The support of the spectrum is $\Supp\widehat f=\{\lambda\in\Lambda_k^n: \widehat f(\lambda)\neq 0\}$.
	
	\begin{Question}\label{qu:sparse} Given a set of $\OO(n)$-types $\Lambda\subset\Lambda_k^n$, and a non-zero distribution $f$ whose spectrum is supported on $\Lambda$, how small can the support of $f$ be?
	\end{Question}
We will provide an obstruction on $\Lambda$ for $f$ to be supported at a point.

\medskip

Let us provide context for this question. The relationship between the zero set sizes of a function and its spectrum is one of several flavors of uncertainty principles in harmonic analysis. For a thorough exposition of this broad area, see \cite{folland_sitaram}, or \cite{wigderson2} for a more recent overview. A typical support-type uncertainty principle asserts that the supports, or possibly the complements of the zero sets, of a function and its spectrum cannot both be too small. In variations on this principle, the support could be replaced by an approximate support. In another type of variation, the rate of decay of the function or of its spectrum may be bounded from below. 

A classical example is the Paley-Wiener theorem, which implies in particular that the Fourier transform of a compactly supported function must have full support. A prototypical theorem due to Benedicks \cite{benedicks_finite_measure} asserts that an integrable function on $\R^n$ and its Fourier transform cannot both have zero sets of finite measure. The entropy uncertainty principle of Beckner \cite{beckner} can be considered as a variation of a support-type uncertainty principle, as the differential entropy of a function can be thought of as the size of its effective support.

In the realm of finite groups, support size is perhaps the most natural measure of localization. Such an uncertainty principle in the setting of finite abelian groups was established by Donoho-Stark \cite{donoho_stark}, see also Tao \cite{tao05}. It was  followed by developments in the non-abelian setting by Meshulam \cite{meshulam06} and others. 

On various non-compact Lie groups, extensions of the theorem of Benedicks exist, see \cite[section 7]{folland_sitaram} and references therein.
In compact settings where the spectrum is discrete, a finitely supported spectrum typically results in an analytic function with full support. Thus interesting generalizations of Benedicks' theorem must consider and quantify infinite spectra, e.g. using their density.

On the circle, such a support-type uncertainty principle is due to Zygmund \cite{zygmund_trigonometric}. It states that a square-integrable function vanishing in a set of positive measure cannot have lacunary Fourier series. A related interesting result is due to Mandelbrojt-Nazarov \cite{mandelbrojt_nazarov}. 
On general compact Lie groups, a recent extension of Zygmund's theorem of Narayanan-Sitaram \cite{narayanan_sitaram_lacunary} asserts that a square-integrable, properly supported function cannot have a lacunary spectrum in a certain sense.

A theorem of Frostman for functions on $S^1$ extended by Beurling \cite{beurling} to distributions, which is however not of pure support-type, asserts that if a distribution is supported on a set of small Hausdorff dimension, its Fourier coefficients cannot decay too fast.  Let us also mention that some results relating the size of the support of a function on $\R$ with the density of the support of its spectrum appeared recently in \cite{amit_olevskii,nazarov_olevskii}.

\medskip

Let us now introduce a prototypical class of functions on the grassmannian that we will consider, which is of geometric origin. Denote by $|\cos(E, E')|$ the cosine of the angle between linear subspaces $E$ and $E'$. It can be defined by $|\cos(E, E')|=\sqrt{\det(P_E\circ P_{E'}:E\to E)}$, where $P_F$ is the orthogonal projection to $F$. 

The cosine transform $\Cos:L^2(\Gr_k(\R^n))\to L^2(\Gr_k(\R^n))$ is the bounded operator $$\Cos(h)(E')=\int_{\Gr_k(\R^n)}|\cos(E, E')| h(E)dE.$$
It is self-adjoint and preserves the class of smooth functions, and so extends to the space of distributions, still denoted $\Cos:C^{-\infty}(\Gr_k(\R^n))\to C^{-\infty}(\Gr_k(\R^n))$.

The cosine transform appears naturally in convex and stochastic geometry, see e.g. \cite{groemer,goodey_howard1, matheron_book, schenider_cosine}. In particular, it plays an important role in integral geometry, namely in convex valuation theory as we explain in section \ref{sec:background}. Furthermore, it fits naturally into the family of $\alpha$-cosine transforms which includes the Radon transform, a central player in integral geometry.
	
The cosine transform is $\OO(n)$-equivariant, and in fact can be rewritten as a $\GL_n(\R)$-equivariant operator between sections of certain line bundles. It is thus only natural that representation theory has been very successfully applied to study it, e.g. in \cite{alesker_bernstein, rubin_cosine, olafsson_cosine, zhang_cosine}. 

In particular, the $\OO(n)$-types appearing in $\textrm{Image}(\Cos)$ have been described in \cite{alesker_bernstein}. At the same time, it is far from clear what are the geometric manifestations of this representation-theoretic description.

\begin{Question}\label{qu:cosine} Given a function, or more generally a distribution $h$ on $\Gr_k(\R^n)$, what are the a-priori restrictions on the support of $\Cos(h)$? 
\end{Question}
We will deduce one such restriction from the quasianalytic property that we will establish for $\textrm{Image}(\Cos)$.

	\subsection{Main results}

	The set $\Lambda_k^n$ of $\OO(n)$-types in $L^2(\Gr_k(\R^n))$ can be identified (see section \ref{sec:background}) with the set $\Lambda_\kappa$ of partitions $\{\lambda_1\geq\dots\geq \lambda_\kappa\}\cap (2\mathbb Z_+)^\kappa$, where $\kappa=\min(k,n-k)$ throughout the paper. Write $|\lambda|=\sum_{i=1}^\kappa \lambda_i$.

	We will need to introduce some terminology.
	
	\begin{Definition}
		The set $\Lambda\subset\Lambda_\kappa$ is \emph{sparse} if $$\lim_{m\to\infty}\frac{|\{\lambda\in \Lambda:|\lambda|\leq2m\}|}{|\{\lambda\in \Lambda_\kappa:|\lambda|\leq2m\}|}=0.$$
		The complement of a sparse set is \emph{co-sparse}.
	\end{Definition}
	Our first main result addresses Question \ref{qu:sparse} with the smallest support imaginable.
	\begin{MainTheorem}\label{thm:sparse}
		Assume $1\leq k\leq n-1$, and $f_0\in C^{-\infty}(\Gr_k(\R^n))$ is supported on a single $E_0\in \Gr_k(\R^n)$. Denote $\Lambda=\Supp\widehat {f_0}\subset \Lambda_{\kappa}$. Then $\Lambda$ is co-sparse. Moreover for $k\in\{1,n-1\}$, $\Lambda$ must be co-finite.
	\end{MainTheorem}
	
	For $2\leq k\leq n-2$, examples of such $f_0$ exist where $\Lambda$ is not co-finite. In section \ref{sec:discussion} we construct such an example which is moreover in the kernel of the cosine transform.
	
	Adapting terminology of harmonic analysis to distributions, one would say that $\{E_0\}\subset \Gr_k(\R^n)$ forms a weakly annihilating pair with any subset of $\Lambda_\kappa$ which is not co-sparse. To our knowledge, support-type uncertainty principles for distributions of lower dimensional support have not been considered in the non-abelian compact setting. Theorem \ref{thm:sparse} can be seen as analogous to Zygmund's theorem, but operating on a different scale.
	\medskip
	Consider now Question \ref{qu:cosine}. For $k\in\{ 1, n-1\},$ the cosine transform is invertible \cite{groemer}, and so there are no restrictions on the support of $\Cos(h)$. For $2\leq k\leq n-2$, the situation is starkly different. We will need
	\begin{Definition}
	For a subspace $F\in\Gr_{n-k}(\R^n)$, define $$\Sigma^k_F=\{E\in\Gr_k(\R^n): E\cap F=\{0\}\}\subset\Gr_k(\R^n).$$ 
	Its complement is $$\Xi^k_F=\Gr_k(\R^n)\setminus \Sigma^k_F=\{E\in\Gr_k(\R^n): E\cap F\neq \{0\}\}.$$
	We write $\Sigma^k,\Xi^k$ if the choice of $F$ plays no role.
	\end{Definition}
	Thus $\Sigma^k_F$ is the unique open Schubert cell with respect to any full flag in $\R^n$ containing $F$. Its complement $\Xi^k_F$ has codimension $1$, and consists of a finite union of locally closed submanifolds.
	
	Given a function $f\in C^\infty(\Gr_k(\R^n))$, we say that $f$ \emph{vanishes exponentially} at $\Xi^k_F$ if $\log|f(E)|\leq b-\frac{C}{d_P(E, \Xi^k_F)}$ for some $C>0$ and $b\in\R$, where $d_P$ is the distance function on $\Gr_k(\R^n)$ induced by any Euclidean structure $P$ on $\R^n$. 
	
	Our second main result is one step towards answering Question \ref{qu:cosine}.
	\begin{MainTheorem}\label{thm:cosine}
		Assume $2\leq k\leq n-2$, and $h\in C^{-\infty}(\Gr_k(\R^n))$. If $\Supp\Cos (h)\subset\Sigma^k$ then $\Cos(h)=0$.
		Moreover, if $h\in C^{\infty}(\Gr_k(\R^n))$ and $\Cos (h)$ vanishes exponentially at $\Xi^k$, then $\Cos(h)=0$. 
	\end{MainTheorem}

It might not be immediately obvious that properly supported functions in the image of $\Cos$ exist at all. Some simple examples of smooth functions $h$ for which $\Supp \Cos (h)$ is a proper subset of arbitrarily small measure can be constructed by exploiting the fact that any $\SO(n-1)$-invariant function is in the image of the cosine transform, which follows from \cite{alesker_bernstein}. More examples can be constructed out of certain $\OO(p,q)$-invariant convex valuations, using descriptions obtained in \cite{alesker_faifman_lorentz, bernig_faifman_opq}. For more details on these examples, see section \ref{sec:discussion}. In all of them, the support is not contractible.
			
The quasianalytic (in the sense of Bernstein for generalized functions, or in the stronger sense of exponential vanishing for smooth sections) property appearing in Theorem \ref{thm:cosine} is shared by several spaces of functions, as we describe below, and thus deserves a name.
\begin{Definition}
	A subspace of generalized functions $\mathcal A\subset C^{-\infty}(\Gr_k(\R^n))$ is \emph{Bernstein $\Xi$-quasianalytic} if, whenever $f\in\mathcal A$ vanishes in a neighborhood of $\Xi^k$, $f=0$.
\\
	A subspace of smooth functions $\mathcal A\subset C^{\infty}(\Gr_k(\R^n))$ is \emph{exponentially $\Xi$-quasianalytic} if, whenever $f\in\mathcal A$ vanishes exponentially at $\Xi^k$, $f=0$.
\\	
	 Both definitions extend immediately to subspaces of generalized, resp. smooth sections of a trivializable line bundle $L$ over $\Gr_k(\R^n)$.
\end{Definition}
Thus Theorem \ref{thm:cosine} asserts that $\Cos(C^{-\infty}(\Gr_k(\R^n)))$ is a Bernstein $\Xi$-quasianalytic class for $2\leq k\leq n-2$, while $\mathrm{Image}\Cos\cap C^{\infty}(\Gr_k(\R^n))$ is exponentially $\Xi$-quasianalytic.

We remark that Theorem \ref{thm:cosine} is false for the kernel of the cosine transform: there are elements in $\Ker(\Cos)$ supported in any open set. This is explained in section \ref{sec:discussion}.

Recall that Aleksandrov's projection theorem \cite{aleksandrov} asserts that a centrally-symmetric convex body in $\R^n$ is uniquely determined by the $k$-volumes of its orthogonal projections to all $k$-dimensional subspaces, for any $1\leq k\leq n-1$. An immediate corollary of Theorem \ref{thm:cosine} in geometric tomography is a sharpening of Aleksandrov's theorem. 
\begin{Corollary}\label{cor:alexandrov}
	Fix $2\leq k\leq n-2$, and let $U\subset\Gr_k(\R^n)$ be a neighborhood of $\Xi^k$. A centrally-symmetric convex body $K\subset\R^n$ of dimension at least $k+1$ is then uniquely determined by its projection function $f_K(E)=\vol_k(P_E(K))$, $E\in U$.
\end{Corollary}
 Previous works on sharpening Aleksandrov's theorem, concerning projections to a subset of hyperplanes, appeared in \cite{grinberg_quinto, schneider_polytopes_projections, schneider_weil_polytopes, schneider_brightness}. Extension to the non-symmetric setting appeared e.g. in \cite{goodey_schneider_weil_determination_projection, goodey_schneider_weil_directed}.

We can also apply Theorem \ref{thm:cosine} to convex valuation theory. Recall that the space of translation-invariant valuations $\Val(\R^n)$ consists of finitely-additive, translation-invariant measures on compact convex sets which are moreover continuous with respect to the Hausdorff metric. Let $\Val_k^+(\R^n)$ denote the subspace of even, $k$-homogeneous valuations. 

Klain's injectivity theorem \cite{klain_even} asserts that $\phi\in \Val_k^+(\R^n)$ is uniquely determined by its restrictions to all $k$-dimensional linear subspaces. Theorem \ref{thm:cosine} implies that one can disregard any compact subset of subspaces inside $\Sigma^k$, as follows. 

\begin{MainTheorem}\label{thm:klain_sharper}
	Assume $2\leq k\leq n-2$, and let $U\subset\Gr_k(\R^n)$ be any open neighborhood of $\Xi^k$. 
	A valuation $\phi\in\Val_k^+(\R^n)$ is uniquely determined by its restrictions to all subspaces $E\in U$.
\end{MainTheorem}
We also prove a similar result for odd valuations, sharpening Schneider's injectivity theorem. A sharper form for odd valuations is conjectured in Section \ref{sec:discussion}.

Combining results of \cite{howe_lee_gln} on certain representations of $\GL_n(\R)$ with more recent results of \cite{alesker_cosine, alesker_gourevitch_sahi} providing explicit descriptions of certain $\GL_n(\R)$-intertwining operators allows us to treat the general $\alpha$-cosine transform $S_\alpha$, which is recalled in subsection \ref{sub:grassmannian}. Theorem \ref{thm:cosine} is but the special case $\alpha=1$ of the following.
\begin{MainTheorem}\label{thm:alpha_cosine}
	Assume $2\leq k\leq n-2$. 
	\begin{enumerate}
		\item If $\alpha\in\mathbb Z$ and $\alpha\geq -\min(k,n-k)+1$, then $\Supp S_\alpha h\subset \Sigma^k$ implies $S_\alpha h=0$, for any $h\in C^{-\infty}(\Gr_k(\R^n))$.  Furthermore, $\mathrm{Image}S_\alpha\cap C^\infty(\Gr_k(\R^n))$ is exponentially $\Xi$-quasianalytic.
		\item  For any other $\alpha\in\mathbb C$ and open subset $U\subset\Gr_k(\R^n)$, there is $h\in C^\infty(\Gr_k(\R^n))$ such that $S_\alpha h\neq 0$ and $\Supp S_\alpha h\subset U$. 
	\end{enumerate}
\end{MainTheorem}

Another related result we obtain is a support theorem for the Radon transform between grassmannians of different dimension, see subsection \ref{sub:grassmannian} for its definition.
\begin{MainTheorem}\label{thm:radon}
	Let $\mathcal R_{p,k}: C^{-\infty}(\Gr_p(\R^n))\to C^{-\infty}(\Gr_k(\R^n))$ be the Radon transform, and $\dim \Gr_p(\R^n)<\dim\Gr_k(\R^n)$. Then $\mathcal R_{p,k}(C^{-\infty}(\Gr_k(\R^n)))$ is Bernstein $\Xi$-quasianalytic, while $\mathrm{Image}\mathcal R_{p,k}\cap C^{\infty}(\Gr_k(\R^n))$ is exponentially $\Xi$-quasianalytic. 
	\end{MainTheorem}
In fact, the same results hold for various $\GL_n(\R)$ modules that appear as subspaces of the space of generalized/smooth sections of certain line bundles over higher rank grassmannians, and form small subspaces in a certain sense. See Remark \ref{rmk:gln_modules} for a precise statement.
	
An immediate corollary in geometric tomography is a sharpening of the Funk section theorem, see Corollary \ref{cor:funk_section}. Previous work extending Funk's theorem to smaller sets of hyperplane sections appeared in \cite{grinberg_quinto}, and in \cite{nazarov_ryabogin_zvavitch} for the non-symmetric case.
	
	\subsection{Overview of the proofs.}
	A counterexample to Theorem \ref{thm:sparse} implies the existence of a differential operator that contains  in its kernel all zonal harmonics in the complement of $\Lambda$. After some manipulations, this becomes a non-trivial PDE with polynomial coefficients, which is solved by all generalized Jacobi polynomials corresponding to zonal harmonics inside the complement of $\Lambda$. It follows that the space of polynomial solutions to this PDE is too large, contradicting Theorem \ref{thm:PDE_main}. 	
	
	The proof of Theorems \ref{thm:cosine}, \ref{thm:alpha_cosine} and \ref{thm:radon} proceeds as follows. Assume a function $f$ in the image of the corresponding integral transform is supported in $\Sigma^k$. First, we use the $\GL_n(\R)$-equivariance of the transform, considered as an intertwining operator between sections of certain line bundles, to shrink the support of $f$ to a point and obtain in the limit, after a carefully chosen rescaling, a distribution $f_0$ supported at a point, which still lies in the image of the transform. Working carefully, one can relax the assumption on the support of $f$ with the weaker assumption that $f$ vanishes at $\Xi^k$ with all derivatives.
	
	We then make use of the description of the $\OO(n)$-types appearing in the image, which is due to Alesker-Bernstein \cite{alesker_bernstein} for the cosine transform, and reduces to results of Alesker et al \cite{alesker_cosine, alesker_gourevitch_sahi} and Howe-Lee \cite{howe_lee_gln} for the $\alpha$-cosine transform, to verify that $\Supp\widehat {f_0}$ is sparse. This then contradicts Theorem \ref{thm:sparse}.
	The proof of the second part of Theorem \ref{thm:alpha_cosine} is by explicit construction of distributions supported at a point, using known representation-theoretic and analytic descriptions of the $\alpha$-cosine transform \cite{alesker_cosine, alesker_gourevitch_sahi, gourevitch_gln}, as well as of the Radon transform \cite{gonzalez_kakehi, kakehi99}.
	
	The proof of Theorem \ref{thm:PDE_main} was explained to us by Joseph Bernstein. A weaker version of this theorem, yielding a correspondingly weaker form of Theorem \ref{thm:sparse} but which is nevertheless sufficient for all the applications in this note, appears in appendix \ref{app:PDE}. Its proof reduces to elementary linear algebra.
		
	\subsection{Plan of the paper}
	
	In section \ref{sec:background} we provide the necessary background from the geometry of grassmannians and the representation theory we will need. In particular, we verify that the $\OO(n)$-types appearing in the images of various integral operators form sparse sets. We also recall some convex valuation theory. In section \ref{sec:sparse} we prove Theorem \ref{thm:sparse}, building on Theorem \ref{thm:PDE_main}. Subsequently in section \ref{sec:cosine} we formulate Theorem \ref{thm:unified}, which is the general statement of which Theorem \ref{thm:radon} and the first part of Theorem \ref{thm:alpha_cosine} are special cases, and prove it by reduction to Theorem \ref{thm:sparse}. We also carry out the verification of the second part of Theorem \ref{thm:alpha_cosine}, and deduce Theorem \ref{thm:klain_sharper}. Then in section \ref{sec:PDE}, which is independent of the rest of the paper, we prove Theorem \ref{thm:PDE_main}. Finally in section \ref{sec:discussion} we provide explicit examples of properly supported functions in the image of the cosine transform originating in convex valuations, as well as examples of distributions supported at a point such that their spectra are not co-finitely supported, in particular one such example inside $\Ker(\Cos)$. We then discuss some further questions and conjectures.

\subsection{Acknowledgements}
Special thanks are due to Joseph Bernstein for explaining to me the proof of Theorem \ref{thm:PDE_main}, as well as to the anonymous referee who identified a significant gap in a previous version of Proposition \ref{prop:point_rescale}, and whose careful reading and numerous suggestions greatly contributed to a better exposition. I am also grateful to Semyon Alesker for bringing reference \cite{gourevitch_gln} to my attention and for valuable comments on the first draft, and to Misha Sodin for some illuminating explanations and references on lacunary Fourier series. The inquisitive questions and remarks of Gil Solanes, and Jan Kotrbaty's comments on a previous version both contributed to the improvement of the paper and are warmly appreciated.

  	\section{Background}\label{sec:background}
  	We write $V$ for a real $n$-dimensional space, and $\R^n$ when the standard Euclidean structure is needed. Let $P_E$ denote the orthogonal projection onto $E$.
  	 
  	\subsection{Grassmannians, line bundles and integral operators}\label{sub:grassmannian}
  	 Recall that $\kappa=\min(k, n-k)$. The principal angles $0\leq \theta_1,\dots,\theta_\kappa\leq \frac\pi2 $ between two subspaces $E,E'\in \Gr_k(\R^n)$ are defined by letting $\cos^2\theta_j$ be the non-trivial eigenvalues of the map $P_E\circ P_{E'}:E\to E$. In particular, $|\cos(E, E')|=\prod_{i=1}^\kappa\cos\theta_j$. For $F=(E')^\perp\in\Gr_{n-k}(\R^n)$, the principal angles between $E$ and $F$ are $\frac \pi 2-\theta_1,\dots, \frac\pi 2-\theta_\kappa$.
  		
  	\begin{Definition}
  		The Radon transform $\mathcal R_{p,k}: C^{\infty}(\Gr_p(\R^n))\to C^{\infty}(\Gr_k(\R^n))$ is given by 
  		
  		$$\mathcal R_{p,k}f (F)=\left\{	\begin{array}{cc}\int_{E\subset F}f(E)dE&,\quad k\geq p\\
  			\int_{E\supset F}f(E)dE&,\quad k\leq p\end{array}\right.$$
  	\end{Definition}
  	It extends by self-adjointness to the space of distributions.
  	\begin{Theorem}[Gelfand-Graev-Rosu  \cite{gelfand_graev_rosu}]\label{thm:radon_injective}
  		When $\dim \Gr_p(\R^n)\leq \dim \Gr_k(\R^n)$, the Radon transform is injective. It is surjective if $\dim \Gr_p(\R^n)\geq\dim \Gr_k(\R^n)$.
  	\end{Theorem}
  	
  	An analytic description of its range is also available. The following statement appears in \cite{kakehi99}.   	
  	\begin{Theorem}[Grinberg \cite{grinberg_radon}, Gonzalez-Kakehi \cite{gonzalez_kakehi}] \label{thm:radon_range_pde}
  		There exists an $\OO(n)$-invariant differential operator $\Omega_{p,k}$ on $\Gr_k(\R^n)$ such that $\Ker(\Omega_{p,k})=\mathrm{Image}(\mathcal R_{p,k})$. 
  	\end{Theorem}
   	
  	The following definitions and facts can be found in \cite{alesker_gourevitch_sahi}.
  	The $\alpha$-cosine transform  $T_\alpha:C^\infty(\Gr_k(\R^n))\to C^\infty(\Gr_k(\R^n))$ is defined for $\alpha\in\mathbb C$, $\mathrm{Re}\alpha>-1$ by
  	$$ T_\alpha(h)(E')=\int_{\Gr_k(\R^n)}|\cos(E, E')|^\alpha h(E) dE,$$
  	where $dE$ is the invariant probability measure.
  	
  	It is well known that the family of operators $T_\alpha$ admits a meromorphic extension in $\alpha\in\C$, see e.g. \cite{rubin_cosine} for an analytic argument, or \cite[Theorem 10.1.6]{wallach} for a general statement concerning meromorphic families of intertwining operators.
  	  	\begin{Definition}
  		The $\alpha_0$-cosine transform $S_{\alpha_0}: C^\infty(\Gr_k(\R^n))\to C^\infty(\Gr_k(\R^n))$ is the leading Laurent coefficient at $\alpha=\alpha_0$:
  		$$T_\alpha=\frac{S_{\alpha_0}}{(\alpha-\alpha_0)^N}+\dots.$$ 
  	\end{Definition}
	In particular, $S_\alpha=T_\alpha$ when $\mathrm{Re}\alpha>-1$, and the cosine transform is $\Cos=S_1=T_1$. Furthermore, $S_\alpha$ has closed image, see e.g. \cite[Corollary 1.6]{gourevitch_gln}. 
		
  	The $\alpha$-cosine transform is symmetric, and extends to an operator $S_\alpha: C^{-\infty}(\Gr_k(\R^n))\to C^{-\infty}(\Gr_k(\R^n))$ with closed image (see e.g. \cite[Claim 4.3]{alesker_faifman_lorentz} for the last claim). We will use the same notation for a given operator with the domain being either smooth or generalized functions; the case at hand should be understood from context. It follows from the above that $$S_\alpha(C^{-\infty}(\Gr_k(\R^n)))\cap C^\infty(\Gr_k(\R^n))=S_\alpha(C^{\infty}(\Gr_k(\R^n))).$$
 
 	It is clear that $S_\alpha$ is $\OO(n)$-equivariant. We will make use of the following result.
 	\begin{Theorem}[Alesker-Gourevitch-Sahi \cite{alesker_gourevitch_sahi}]\label{thm:cosine_diff_operator}
 		For any $\alpha\notin [-(\kappa+1),-2]\cap \mathbb Z$ there is an $\OO(n)$-invariant differential operator $D_{\alpha,\alpha+2}$ such that
 		$S_\alpha=D_{\alpha,\alpha+2} \circ S_{\alpha+2}$.\\
 		For $\alpha\in[-\kappa, -1]\cap \mathbb Z$, $S_\alpha$ is the composition of two Radon transforms with the intermediate grassmannian being $\Gr_{|\alpha|}(\R^n)$.
 	\end{Theorem}	
 
 	 In \cite{alesker_bernstein, alesker_cosine} it was observed that $S_\alpha$ can be made $\GL_n(\R)$-equivariant by considering it as on operator between sections of certain line bundles, which we now recall. 
 	
 	Denote the line of $\alpha$-densities on a linear space $E$ by $\Dens^\alpha(E)$. 
 	Let $M_\alpha, L_\alpha$ be two line bundles over $\Gr_k(V)$, with fibers given by
 	$$M_\alpha|_E=\Dens^\alpha(E),\quad L_\alpha|_E=\Dens^\alpha(V/E)\otimes \Dens(T_E\Gr_k(V)).$$
 	We will make use of the equivalent description
 	\begin{equation}\label{eq:L_alpha} L_\alpha|_E=\Dens^{-n-\alpha}(E)\otimes \Dens^{k+\alpha}(V). \end{equation}
 	
 	We will write $\Gamma^\infty(\Gr_k(V), L)$ for the smooth sections of the line bundle $L$, and $\Gamma^{-\infty}(\Gr_k(V), L)$ for the generalized sections.
 	
 	Observe that $L_\alpha|_E=M_\alpha|_E^*\otimes \Dens(T_E\Gr_k(V)) \otimes\Dens^\alpha(V)$ and so there is a non-degenerate pairing 
 	$$\Gamma^\infty(\Gr_k(V), M_\alpha)\times \Gamma^\infty(\Gr_k(V), L_\alpha\otimes \Dens^{-\alpha}(V))\to \R$$
 	  	given by integrating the product over $\Gr_k(V)$.
 	  	
 	An important special case is $\alpha=0$. Then $\Gamma^\infty(\Gr_k(V), L_0)=\mathcal M^\infty(\Gr_k(V))$ is the space of smooth measures on $\Gr_k(V)$, while $\Gamma^{-\infty}(\Gr_k(V), M_0)=C^{-\infty}(\Gr_k(V))$ is its topological dual, the generalized functions on $\Gr_k(V)$.  	
 	  	
 	For $g\in \GL(V)$ and $\phi\in   \Gamma^\infty(\Gr_k(V), M_\alpha)$, the action of $g$ on $\phi$ is denoted $g\phi$, and given by $g\phi (E)= g_*(\phi(g^{-1}E))$. The action on $\mu\in \Gamma^\infty(\Gr_k(V), L_\alpha\otimes \Dens^{-\alpha}(V))$ is by duality, namely $\langle \phi, g\mu\rangle =\langle g^{-1}\phi, \mu\rangle$.
 	  	
 	The $\alpha$-cosine transform can be rewritten as a $\GL(V)$-equivariant operator $$S_\alpha:\Gamma^\infty(\Gr_k(V), L_\alpha)\to\Gamma^\infty(\Gr_{n-k}(V), M_\alpha),$$ which coincides with the previous definition once a Euclidean structure is used to trivialize both line bundles, and to identify $\Gr_k(V)$ with $\Gr_{n-k}(V)$ using the orthogonal complement map.
   
 As with the $\alpha$-cosine transform, $\mathcal R_{p,k}$ can be rewritten  $\GL_n(\R)$-equivariantly \cite{gelfand_graev_rosu}. E.g. for $k\geq p$ one obtains
 $$\mathcal R_{p,k}: \Gamma^{\pm \infty}(\Gr_p(V), M_{-k})\to\Gamma^{\pm \infty}(\Gr_k(V), M_{-p}).$$

  	\subsection{Polynomials and representation theory}
  	Let $\mathcal P_m\subset \mathbb C[x_1,\dots, x_k]$ denote the subspace of polynomials of degree at most $m$.
  
  	We first collect some basic facts on polynomials.
  	  	
  	\begin{Lemma} \label{lem:polynomials_dimension} 
 	 	It holds that $\dim \mathcal P_m = \frac{1}{k!}m^k+ O(m^{k-1})$ as $m\to\infty$.
   	\end{Lemma}
  \proof
  	Denote by $\mathcal H_m\subset \mathcal P_m$ the subspace of homogeneous polynomials of degree $m$. It holds that 
  $\mathcal P_m=\mathcal H_0\oplus \mathcal H_1\oplus \dots\oplus \mathcal H_m$.
  
   It holds that $\dim\mathcal H_m={m+k-1\choose k-1}=\frac{(m+1)\cdots(m+k-1)}{(k-1)!}$, in particular $$\frac{1}{(k-1)!}m^{k-1}\leq \dim\mathcal H_m\leq \frac{1}{(k-1)!}m^{k-1} +O(m^{k-2}).$$ The stated asymptotics is now immediate to deduce.
  \endproof
  We will frequently make use of the subset $\Lambda_k(2m):=\{\lambda\in \Lambda_k: |\lambda|\leq 2m\}$.
  Denote by $\mathcal P_m^s[y_1,\dots,y_k]\subset\mathbb C[y_1,\dots,y_k]$ the symmetric polynomials of degree at most $m$. 
  \begin{Lemma}\label{lem:symmetric_dimension}
	 It holds that
	 \begin{equation}\label{eq:symmetric_asymptotics}|\Lambda_k(2m)|=\dim \mathcal P_m^s[y_1,\dots,y_k]= \frac{1}{k!^2}m^k+O(m^{k-1}),\quad m\to\infty.\end{equation}
  \end{Lemma}
\proof
The monomial symmetric polynomials $$\left\{\sum_{\sigma\in S_k} y_{\sigma 1}^{\lambda_1/2}\cdots  y_{\sigma k}^{\lambda_k/2}\right\}_{\lambda\in\Lambda_k(2m)}$$ form a basis of $\mathcal P_m^s[y_1,\dots,y_k]$, implying the first equality.

The asymptotics follow from standard facts on the restricted partition function. Denoting $P(m,k)=|\{\lambda\in\Lambda_k: |\lambda|=2m\}|$, it holds that 
$$\frac{1}{k!}{m+k-1 \choose k-1}\leq P(m,k)\leq \frac{1}{k!}{m+{k + 1 \choose 2}-1 \choose k-1},$$
see e.g. \cite{szekeres_partitions}. Thus $P(m,k)=\frac{1}{k!(k-1)!}m^{k-1}+O(m^{k-2})$ as $m\to\infty$. 

Since $\sum_{j=1}^m j^{k-1}=\frac{m^k}{k}+O(m^{k-1})$ as $m\to\infty$, we deduce
\[|\Lambda_k(2m)|=P(0,k)+P(1,k)+\dots+P(m,k)=\frac{1}{k!^2}m^k+O(m^{k-1}), \]
as claimed.
\endproof
   We next discuss representations of $\OO(n)$. The following can be found in \cite{sugiura, zhelobenko}. 
   
    Recall that $\kappa=\min(k,n-k)$. The real grassmannian is a symmetric space, and consequently we can write $$L^2(\Gr_k(\R^n))=\oplus_{\lambda\in \Lambda_\kappa} V_\lambda,$$ where the $V_\lambda$ are pairwise non-isomorphic, finite-dimensional irreducible representations of $\OO(n)$. 
     	
  	Fixing $E_0\in\Gr_k(\R^n)$ with its stabilizer $H=\OO(k)\times\OO(n-k)$, each $V_\lambda$ contains a unique one-dimensional subspace which is $H$-invariant. It is spanned by the zonal harmonic $Z_\lambda$, which has the property that $\langle f, Z_\lambda\rangle = \widehat f(\lambda)(E_0)$ for any $f\in C^\infty(\Gr_k(\R^n))$. 
  	
  	The zonal harmonics $Z_\lambda\in V_\lambda$ are given by $Z_\lambda(E)=P_\lambda(y_1,\dots, y_\kappa)$, where $y_j=\cos^2\theta_j(E, E_0)$, and $P_\lambda$ are certain symmetric polynomials called the generalized Jacobi polynomials \cite{james_constantine}. It holds that $\deg P_\lambda=\frac12|\lambda|$.

   	\begin{Lemma}\label{lem:symmetric_basis}
  		The set $\{P_\lambda(y_1,\dots,y_\kappa)\}_{\lambda\in\Lambda_\kappa(2m)}$ is a basis of $\mathcal P_m^s[y_1,\dots,y_\kappa]$.
  	\end{Lemma}
  \proof
  	As $P_\lambda\in \mathcal P_m^s[y_1,\dots,y_\kappa]$ when $|\lambda|\leq 2m$, and all $P_\lambda$ are linearly independent, this follows at once from Lemma \ref{lem:symmetric_dimension}.
  \endproof
  	
  Let us record also the following trivial fact.
  \begin{Lemma}\label{lem:finite_support}
  	If $f,g\in C^\infty(\Gr_k(\R^n))$ both have finitely supported spectrum, then so does $f g$.
  \end{Lemma}
\proof
Let $W_f, W_g\subset  C^\infty(\Gr_k(\R^n))$ be finte-dimensional $\OO(n)$-invariant subspaces containing $f$ and $g$, respectively. 
Choose bases $f_1,\dots, f_l$ for $W_f$, resp. $g_1,\dots, g_m$ for $W_g$.
Then $fg\in \Span\{f_ig_j: 1\leq i\leq l, 1\leq j\leq m\}$, and the latter space is evidently $\OO(n)$-invariant and finite-dimensional.
\endproof

\begin{Proposition}\label{prop:radon_types}
	Assume $\dim\Gr_p(\R^n)<\dim \Gr_k(\R^n)$. The $\OO(n)$-types in $\mathrm{Image}(\mathcal R_{p,k})$ are $\{\lambda\in\Lambda_\kappa: \lambda_{\min(p,n-p)+1}=0\}$.
\end{Proposition}
\proof
This is an immediate consequence of the injectivity or the Radon transform \cite{gelfand_graev_rosu} and the description of the $\OO(n)$-types in $L^2(\Gr_k(\R^n))$.
\endproof

\begin{Lemma}\label{lem:radon_image}
	The $\OO(n)$-types appearing in the image of $\mathcal R_{p,k}: C^{\infty}(\Gr_p(\R^n))\to C^{\infty}(\Gr_k(\R^n))$ form a sparse set when 
	$\dim\Gr_p(\R^n)<\dim \Gr_k(\R^n)$.
\end{Lemma}
\proof
We simply ought to check that $\min(p,n-p)<\min(k,n-k)$ implies
$$ \lim_{m\to \infty}\frac {|\Lambda_{\min(p,n-p)}(2m)|}{|\Lambda_{\min(k,n-k)}(2m)|}=0,$$
which follows from Lemma \ref{lem:symmetric_dimension}. \endproof

Next we describe the image of certain $\alpha$-cosine transforms.
\begin{Proposition}[Alesker \cite{alesker_cosine}, Alesker-Gourevitch-Sahi \cite{alesker_gourevitch_sahi}]\label{prop:cosine_types}
It holds that $\mathrm{Image}(S_\alpha)$ consists of the $\OO(n)$-types $\lambda \in A_\alpha$, where 
\begin{itemize}
	\item  for $\alpha\notin \mathbb Z$, $S_\alpha$ is invertible, that is $A_\alpha=\Lambda_\kappa$.
	\item  for $\alpha\in(2\mathbb Z_++1)$, $A_\alpha=\{\lambda\in\Lambda_\kappa: \lambda_2\leq 1+\alpha\}$. 
	\item  for $\alpha\in2\mathbb Z_+$, $A_\alpha=\{\lambda\in\Lambda_\kappa: \lambda_1\leq \alpha\}$.
	\item  for $\alpha\in [-(\kappa-1), -1]\cap \mathbb Z$, $A_\alpha=\{\lambda\in\Lambda_\kappa: \lambda_{|\alpha|+1}=0\}$. 
\end{itemize}	
\end{Proposition}
\proof
The first two items appear in \cite[Theorem 4.15]{alesker_cosine}. The third is easy to deduce from \cite[Theorem 1.6]{gourevitch_gln} and \cite[Theorem 3.4.2]{howe_lee_gln}. The last item follows from the second part of Theorem \ref{thm:cosine_diff_operator}, combined with Theorem \ref{thm:radon_injective}.
\endproof

We now verify that $\mathrm{Image}(S_\alpha)$ consists of a sparse set of $\OO(n)$-types for all integer $\alpha$ in Proposition \ref{prop:cosine_types}.

\begin{Lemma}\label{lem:cosine_image}
	Assume $2\leq k\leq n-2$ and $\alpha\in\mathbb Z$. The $\OO(n)$-types occurring in $\mathrm{Image}(S_\alpha)$ form a sparse set whenever $\alpha\geq-(\kappa-1)$.
\end{Lemma}
\proof
For $\alpha\in2\mathbb Z_+$ the set of types is finite and therefore sparse. \\
For $\alpha\in2\mathbb Z_++1$ we have 
$|\{\lambda\in \Lambda_\kappa(2m): \lambda_2\leq 1+\alpha\}|\leq C_{\alpha, \kappa}(m+1)$, and so this is a sparse set of types by Lemma \ref{lem:symmetric_dimension}.

For $-(\kappa-1)\leq \alpha\leq -1$ this follows from Lemma \ref{lem:radon_image} and Theorem \ref{thm:cosine_diff_operator}.
\endproof
\begin{Remark}
	One could verify directly that for all other values of $\alpha$, the $\OO(n)$-types in the image are co-sparse. This however will follow from the second part of Theorem \ref{thm:alpha_cosine} in light of Theorem \ref{thm:sparse}.
\end{Remark}

	\subsection{Convex valuation theory}
Let $\mathcal K(V)$ be the set of non-empty compact convex sets in $V$, equipped with the topology of the Hausdorff metric.

The translation-invariant continuous valuations $\Val(V)$ consist of translation-invariant, continuous functions $\phi:\mathcal K(V)\to \mathbb C$ satisfying the valuation property
$$\phi(K\cup L)+\phi(K\cap L)=\phi(K)+\phi(L)\quad \forall K,L, K\cup L \in \mathcal K(V).$$ 	
The important subspace of smooth valuations $\Val^\infty(V)$, introduced by Alesker, consists of the smooth vectors for the natural action of $\GL(V)$ on $\Val(V)$, equipped with a certain natural Banach space topology.

By $\Val_k(V)$ we denote the set of $k$-homogeneous valuations. It is a theorem of McMullen \cite{mcmullen_decomposition} that 
$\Val(V)=\oplus_{k=0}^n \Val_k(V).$
The spaces $\Val_0(V)$ and $\Val_n(V)$ are each one-dimensional, and consist of constant valuations and Lebesgue measures, respectively. For $1\leq k\leq n-1$, $\Val_k(V)$ is infinite-dimensional. By $\Val^\pm(V)$ we denote the even (resp. odd) valuations, namely those satisfying $\phi(-K)=\pm \phi(K)$.

Klain's injectivity theorem \cite{klain_even} asserts that the map $$\Kl:\Val_k^+(V)\to \Gamma(\Gr_k(V), M_1),\quad \Kl(\phi)(E)=\phi|_E$$ is injective. Moreover, by \cite[Theorem 1.3]{alesker_bernstein} its image on smooth valuations coincides with the image of the cosine transform on smooth sections.  

Schneider's injectivity theorem \cite{schneider_simple} plays a similar role for odd valuations. It asserts that the map $$\Sc:\Val_k^-(V)\to \Gamma(\Gr_{k+1}(V), \Val_k^-(E)),\quad \Sc(\phi)(E)=\phi|_E$$
is injective.

Klain's and Schneider's theorems have proven remarkably useful in integral geometry. Notably, they were instrumental in Alesker's resolution of McMullen's conjecture \cite{Alesker:Irreducibility}.

The generalized valuations $\Val^{-\infty}(V)$ can be defined as the continuous dual to $\Val^\infty(V)\otimes \Dens(V)^*$, with the latter space equipped with the G\r{a}rding topology. The Alesker-Poincare pairing \cite{alesker_product} guarantees there is a dense embedding $\Val(V)\hookrightarrow\Val^{-\infty}(V)$. The Klain map extends to an injective map  $\Kl:\Val_k^{-\infty,+}(V)\to \Gamma^{-\infty}(\Gr_k(V), M_1)$.

For a comprehensive introduction to convex valuation theory, see \cite{schneider_book, alesker_fu}. For more details on generalized valuations, see e.g. \cite{alesker_faifman_lorentz}.
 	 
\section{Distributions supported at a point}\label{sec:sparse}
Let $\delta_{E_0}\in C^{-\infty}(\Gr_k(\R^n))$ be given by $\langle\delta_{E_0}, \mu\rangle=\frac{d\mu}{dE}(E_0)$ for smooth measures $\mu\in\mathcal M^\infty(\Gr_k(\R^n))$. By $U(\mathfrak g)$ we denote the universal enveloping algebra of the Lie algebra $\mathfrak g$.
 \begin{Lemma}\label{lem:universal}
 	Any non-zero generalized function $f_0\in C^{-\infty}(\Gr_k(\R^n))$ supported at a point $E_0$ is given by $f_0=D\delta_{E_0}$, where $D\in U(\mathfrak {so}_n)$. 
 	\end{Lemma}
\proof
It is well-known \cite{hormander_vol1} that one can write $f_0=D_1\delta_{E_0}$ for some differential operator $D_1$ of order $m$. Let us show by induction on $m$ that we may replace $D_1$ by an element of $U(\mathfrak{so}_n)$ without altering $D_1\delta_{E_0}$. 

For $m=0$, $D_1\delta_{E_0}=\psi(E)\delta_{E_0}$ for some smooth function $\psi$. Since $f_0\neq 0$, $\psi(E_0)\neq 0$ and $f_0=\psi(E_0)\delta_{E_0}$. Thus we may take $D$ to be the constant $\psi(E_0)$. 

Now consider general $m$. Recall that $\langle L_X\delta_{E_0}, \mu\rangle = -\langle \delta_{E_0}, L_X\mu\rangle$ for a smooth measure $\mu$. Therefore writing  $\langle D\delta_{E_0}, \mu\rangle = \langle \delta_{E_0}, D^*\mu\rangle$ for the adjoint differential operator, it holds that $D\in U(\mathfrak{so}_n)$ if and only if $D^*\in U(\mathfrak {so}_n)$. 

 Thus for any smooth measure $\mu$ we have $$\langle D_1\delta_{E_0}, \mu\rangle= \langle \delta_{E_0},  D_1^*\mu\rangle=\frac{d( D_1^*\mu)}{dE}(E_0).$$
Furthermore, one can write $$\frac{d( D_1^*\mu)}{dE}(E_0)=D_1^*(\frac{d\mu}{dE})(E_0)+\frac{d( D_1'\mu)}{dE}(E_0)$$ for some linear differential operator $D_1'$ of order at most $m-1$. 

By the induction hypothesis, it suffices to prove the following statement: for any differential operator $\widehat D_1$ of order $m$, there is $\widehat D\in U(\mathfrak{so}_n)$ such that $\widehat Dh(E_0)=\widehat D_1h(E_0)$  for any $h\in C^\infty(\Gr_k(\R^n))$. We do that again by induction on $m$.

For $m=0$ this is trivial. For general $m$, we may assume by linearity that $\widehat D_1=L_{X_1}\dots L_{X_m}$ for some vector fields $X_j$. 
Find $Y_1\in \mathfrak{so}_n$ such that $Y_1|_{E_0}=X_1|_{E_0}$. Then 
$$\widehat D_1h(E_0)=L_{Y_1}L_{X_2}\dots L_{X_m}h(E_0)=L_{X_2}\dots L_{X_m}L_{Y_1}h (E_0)+ D_1'h(E_0),$$
where $D'_1$ has order at most $m-1$. Applying the induction hypothesis to both $L_{X_2}\dots L_{X_m}$ and $D'_1$ now completes the proof. 

\endproof

\begin{proof}[Proof of Theorem \ref{thm:sparse}] Considered as a map $C^\infty(\Gr_k(\R^n))\to C^\infty(\Gr_{n-k}(\R^n))$, the pull-back by orthogonal complement defines an isomorphism between each pair of irreducible components of $\OO(n)$. We may thus assume that $k\leq \frac n 2$. Denote by $H=\OO(k)\times \OO(n-k)$ the stabilizer of $E_0$.

By Lemma \ref{lem:universal} we may find $D\in U(\mathfrak{so}_n)$ such that $f_0=D\delta_{E_0}$.
The zonal harmonic at $E_0$, denoted $Z_\lambda\in V_\lambda$, is given by $\langle \phi, Z_\lambda\rangle=\widehat \phi(\lambda)(E_0)$. We then have the weakly convergent decomposition $\delta_{E_0}=\sum_{\lambda\in\Lambda_k}  Z_\lambda$.

The operator $D$ leaves all the irreducible components $V_\lambda$ invariant. It therefore holds that $DZ_\lambda=0$ for all $\lambda\in \Lambda^c$, while $DZ_\lambda\neq 0$ for some $\lambda\in \Lambda$.

For any zonal harmonic $Z$ on $\Gr_k(\R^n)$ we can write $Z(E)=f(\sigma_1(E),\dots, \sigma_k(E))$, where $\sigma_j(E)$ are the elementary symmetric polynomials in the squared cosines of the principal angles between $E$ and $E_0$, and $f(\sigma_1,\dots,\sigma_k)$ is a polynomial, namely a generalized Jacobi polynomial written in terms of the elementary symmetric polynomials. Then $$DZ=\sum_{|\alpha|\leq N}\psi_\alpha(E) \partial^\alpha f(\sigma_1(E),\dots,\sigma_k(E)),$$ where $\alpha\in\mathbb Z_+^k$, $|\alpha|=\sum_i \alpha_i$, and each $\psi_\alpha(E)$ is a linear combination of products of derivatives by elements of $U(\mathfrak{so}_n)$ of the various $\sigma_j(E)$. By Lemma \ref{lem:symmetric_basis}, each $\sigma_j(E)$ has finitely supported spectrum. It follows by Lemma \ref{lem:finite_support} that also all $\psi_\alpha$ have finitely supported spectrum. 

We would like to replace $\psi_\alpha$ by their $H$-invariant components, and so we intend to integrate the equation $DZ=0$ over $H$. However, to ensure we do not trivially get $0$ after the averaging, we first choose any $\alpha_0$  appearing in the sum, and multiply the equation $DZ=0$ by $\psi_{\alpha_0}(E)$, obtaining
$$ \sum_{|\alpha|\leq N}\psi_\alpha(E)\cdot \psi_{\alpha_0}(E) \cdot \partial^\alpha f(\sigma_1(E),\dots,\sigma_k(E))=0.$$ 
Again by Lemma \ref{lem:finite_support}, all the functions $\psi_\alpha(E)\cdot \psi_{\alpha_0}(E)\in C^\infty(\Gr_k(\R^n))$ still have finitely supported spectrum.

Now the coefficient of
$\partial^{\alpha_0}f$ is nonzero and non-negative, while all $\sigma_j(E)$ are $H$-invariant, and thus we can integrate over $H$ and get the equation
\begin{equation}\label{eq:impossible_PDE0}  \sum_{|\alpha|\leq N} q_\alpha(E)\partial^\alpha f(\sigma_1(E),\dots,\sigma_k(E))=0,\end{equation}
where $$q_\alpha(E)=\int_H h^*(\psi_\alpha(E)\cdot \psi_{\alpha_0}(E))dh,$$
and the coefficient $q_{\alpha_0}$ of $\partial ^{\alpha_0} f$ remains non-negative and nonzero.

Observe that each $q_\alpha(E)$ is $H$-invariant and has finitely supported spectrum. Thus it is a finite linear combination of zonal harmonics, and we may write $q_\alpha=p_\alpha(\sigma_1(E),\dots,\sigma_k(E))$, where $p_\alpha(\sigma_1,\dots,\sigma_k)$ is a polynomial. 

For $\lambda\in \Lambda_k$, let $\widetilde P_\lambda(\sigma)\in\mathbb R[\sigma_1,\dots,\sigma_k]$ denote the corresponding generalized Jacobi polynomial, written in the basis of elementary symmetic polynomials, so that $Z_\lambda(E)=\widetilde P_\lambda(\sigma_1(E),\dots, \sigma_k(E))\in V_\lambda$. For $\Lambda'\subset\Lambda_k$ write $K_{\Lambda'}[\sigma]=\Span\{\widetilde P_\lambda:\lambda \in \Lambda'\} \subset \mathbb C[\sigma_1,\dots,\sigma_k]$. 

Thus eq. \eqref{eq:impossible_PDE0} is a nonzero PDE in the variables $\sigma_1,\dots,\sigma_k$ with polynomial coefficients:
\begin{equation}\label{eq:impossible_PDE}  \widetilde D(f):=\sum_{|\alpha|\leq N}p_\alpha(\sigma_1,\dots,\sigma_k)\cdot \partial^\alpha f(\sigma_1,\dots,\sigma_k)=0,\end{equation}
and \begin{equation}\label{eq:kernel_subset}K_{\Lambda^c}[\sigma]\subset \Ker\widetilde D.\end{equation}

Now for $k=1$ this is a linear ODE with polynomial coefficients, and so $\Lambda^c$ must be finite.
Assume $2\leq k\leq \frac n 2$.

There is a natural isomorphism of algebras $s: \mathbb C[\sigma_1,\dots, \sigma_k] \to \mathbb C_{\textrm{sym}}[y_1,\dots, y_k]$, where $\sigma_i$ is assigned the $i$-th elementary symmetric polynomial in $y$. It holds that
$ \mathcal P^s_m[y]\subset s(\mathcal P_m[\sigma])$.

By construction, $s(K_{\Lambda'}[\sigma]\cap \mathcal P_m[\sigma])=s(K_{\Lambda'}[\sigma])\cap s(\mathcal P_m[\sigma])$. Recall also that $\deg s(\widetilde P_\lambda)=\frac12|\lambda|$.
 
Assume in contradiction that $\Lambda$ is not co-sparse, that is for some $\epsilon>0$ and $m_j\to\infty$ one has $|\Lambda^c\cap \Lambda_k(2m_j)|\geq \epsilon|\Lambda_k(2m_j)|$.

Recall that by Lemma \ref{lem:symmetric_basis}, $\{s(\widetilde P_\lambda): \lambda\in \Lambda_k(2m)\}$ is a basis of $\mathcal P^s_m[y]$.
Hence \begin{align*}\dim K_{\Lambda^c}[\sigma]\cap \mathcal P_{m_j}[\sigma]&=\dim s(K_{\Lambda^c}[\sigma]\cap \mathcal P_{m_j}[\sigma])=\dim K_{\Lambda^c}[y]\cap s(\mathcal P_{m_j}[\sigma])\\&\geq  
	\dim K_{\Lambda^c}[y]\cap \mathcal P_{m_j}^s[y]=|\Lambda^c\cap \Lambda_k(2m_j)|\\&\geq \epsilon |\Lambda_k(2m_j)|=\epsilon \dim \mathcal P^s_{m_j}[y].\end{align*}
By Lemmas \ref{lem:polynomials_dimension} and \ref{lem:symmetric_dimension}, $\dim \mathcal P_m[\sigma]\sim \frac{1}{k!}m^k$ while $\dim \mathcal P^s_{m}[y]\sim \frac{1}{k!^2}m^k$ as $m\to\infty$. Therefore by eq. \eqref{eq:kernel_subset}, $$\dim \Ker \widetilde D\cap \mathcal P_{m_j}[\sigma]\geq \frac{\epsilon}{k!}(1+o(1)) \dim \mathcal P_{m_j}[\sigma],\quad j\to\infty.$$ But this contradicts Theorem \ref{thm:PDE_main}.
\end{proof}

\section{The support of generalized sections of equivariant line bundles}\label{sec:cosine}
Recall that $V$ is a real $n$-dimensional linear space, and $$\Xi^k_F=\Gr_k(V)\setminus\Sigma^k_F=\{E\in\Gr_k(V): E\cap F\neq \{0\}\}.$$ 

First we observe that $\Sigma^k\subset\Gr_k(V)$ is not only contractible, but contractible by linear transformations. In the next proposition, we will show that given a section of an equivariant line bundle which vanishes exponentially at $\Xi^k$, we can use the action of $\GL(V)$ to shrink $\Sigma^k$ to a point and obtain in the limit,  after a careful rescaling, a nontrivial generalized section that is supported at one point.

Choose any $E_0\in \Sigma^k_F$, and fix a Euclidean structure on $V$ such that $F=E_0^\perp$. For the proof of the following, we will make frequent use of various facts on the geometry of the grassmannian appearing in \cite{kozlov13} (see \cite{kozlov_eng12,kozlov_english_3} for the English version).
\begin{Proposition}\label{prop:point_rescale}
	Let $g_\epsilon\in \End(V)$ be given by $g_\epsilon=P_{E_0}+\epsilon P_{F}$. Assume $\phi$ is a smooth section of $M_\alpha$, and $\phi$ vanishes exponentially at $\Xi^k_F$. Then one can choose coefficients $\lambda_\epsilon\in \R$ such that $\lambda_\epsilon g_\epsilon\phi\to \phi_0$ as $\epsilon\to 0^+$, where $\phi_0$ is a non-zero generalized section of $M_\alpha$ supported at $E_0$.
\end{Proposition}
\proof
 We will use the Euclidean structure to identify $\phi$ with $f\in C^\infty(\Gr_k(V))$, and write $[\phi]=f$. We will use the same convention also for sections of $L_\alpha\otimes\Dens^{-\alpha}(V)$. Denote by $d_P$ the distance function on $\Gr_k(V)$ equipped with the standard Riemannian metric induced by the Euclidean structure. 
 
 It holds that $\det g_\epsilon=\epsilon^{n-k}$.
 Denoting by $\textrm{Jac}$ the Jacobian, consider $\eta_g(E)=\textrm{Jac}(g: E\to gE)^{-1}$, which therefore satisfies $\lim_{\epsilon\to 0}\eta_{g_\epsilon}(E)=|\cos( E, E_0)|^{-1}$.
 
 Denote $f_\epsilon=[g_\epsilon(\phi)]$, and let $u$ be a smooth test function on $\Gr_k(V)$. If $\mu\in\Gamma^\infty(\Gr_k(V), L_\alpha\otimes \Dens^{-\alpha}(V))$ is given by $[\mu]=u$, 
 then by eq. \eqref{eq:L_alpha},
 \begin{align*} \langle f_\epsilon, u\rangle=\int_{\Gr_k(V)} f_\epsilon(E)u(E)dE&=\langle g_\epsilon \phi, \mu\rangle=\langle \phi, g_\epsilon^{-1}\mu\rangle
 	\\&=\int_{\Gr_k(V)} f(E)\eta_{g_\epsilon}(E)^{n+\alpha}\det(g_\epsilon)^ku(g_\epsilon E)dE\\&=\epsilon^{k(n-k)}\int_{\Gr_k(V)} f(E)u(g_\epsilon E)\eta_{g_\epsilon}(E)^{n+\alpha}dE.\end{align*}
 
 Take $\epsilon<1$, and introduce the neighborhood $W_\epsilon$ of $\Xi^k_F$ given by 
 $$W_{\epsilon}=\{E\in\Gr_k(V): d_P(E, \Xi^k_F)\leq \epsilon^{\frac 12} \}.$$ 
 Write 
 
 $$\langle f_\epsilon, u\rangle=\epsilon^{k(n-k)}(I_\epsilon+J_\epsilon),$$
 where 
 
 $$I_\epsilon=\int_{W_\epsilon} f(E)u(g_\epsilon E)\eta_{g_\epsilon}(E)^{n+\alpha}dE,\qquad J_\epsilon=\int_{W_\epsilon^c} f(E)u(g_\epsilon E)\eta_{g_\epsilon}(E)^{n+\alpha}dE.$$
 
 Denote by $0\leq \beta_1,\dots,\beta_\kappa\leq \frac\pi 2$ the principal angles between $E$ and $F$. Then $d_P(E, \Xi^k_F)=\min\beta_i$ \cite{kozlov13}, and $|\cos(E,E_0)|=\prod_{i=1}^\kappa \sin\beta_i$. Since $f$ vanishes exponentially at $\Xi^k_F$, it holds for some $\hat c_0, \hat c_1>0$ that 
 \begin{equation}\label{eq:f_bound}|f(E)|\leq \hat c_0 e^{-\hat c_1/d_P(E, \Xi^k_F)}=\hat c_0 e^{-\hat c_1/\min\beta_i}.\end{equation}
 
 It follows that for all $E\in W_\epsilon$, $|f(E)|\leq \hat c_0 e^{-\hat c_1/\sqrt{\epsilon}}$. Also
 $$   \epsilon^{2\kappa}\leq \prod_{i=1}^\kappa(\sin^2\beta_i+\epsilon^2\cos^2\beta_i)\leq 1,$$
 so that \begin{equation}\label{eq:eta}\eta_{g_\epsilon}(E)=\prod_{i=1}^\kappa(\sin^2\beta_i+\epsilon^2\cos^2\beta_i)^{-1/2}\end{equation}
 satisfies $1\leq \eta_{g_\epsilon}(E)\leq \epsilon^{-\kappa}$.
 
 Thus
 
 \begin{equation}\label{eq:I_epsilon}|I_\epsilon|=| \int_{W_\epsilon}f(E)u(g_\epsilon E)\eta_{g_\epsilon}(E)^{n+\alpha}dE|\leq \hat c_0\|u\|_\infty e^{-\hat c_1/\sqrt{\epsilon}} \cdot \epsilon^{-\kappa|(n+\alpha)|}.\end{equation} 
 
 Let us now examine $J_\epsilon$.  Write 
 $$\eta_{g_\epsilon}(E)^{n+\alpha}=\left(\prod_{i=1}^\kappa \sin^2\beta_i \prod_{i=1}^\kappa(1+\epsilon^2\cot^2\beta_i)\right)^\gamma,\qquad \gamma=-\frac{n+\alpha}{2},$$
 and define $$f_1(E)=f(E)(\prod_{i=1}^\kappa \sin^2\beta_i)^\gamma.$$ 
 Note that \[(\sin d_P(E, \Xi^k_F))^{2\kappa}=(\sin(\min\beta_i))^{2\kappa}\leq \prod_{i=1}^\kappa \sin^2\beta_i\leq 1 ,\] and so by eq. \eqref{eq:f_bound} we have for some $c_0', c_1'>0$ that \begin{equation}\label{eq:f_1_first_bound}|f_1(E)|\leq c_0'e^{-c'_1/d_P(E, \Xi^k_F)}.\end{equation} Define also
 \begin{equation}\label{eq:zeta}\zeta_\epsilon(E)=\prod_{i=1}^\kappa(1+\epsilon^2\cot^2\beta_i)^\gamma,\end{equation}
 so that  \begin{equation}\label{eq:J_epsilon0}J_\epsilon=\int_{W_\epsilon^c} f_1(E)u(g_\epsilon(E))\zeta_\epsilon(E)dE.\end{equation}
 Note that for $E\in W_\epsilon^c$, $\beta_i\geq d_P(E, \Xi_F^k)\geq \sqrt \epsilon$ for all $1\leq i\leq \kappa$, and so for all $i$, \begin{equation}\label{eq:cot_beta}\cot\beta_i\leq \frac{1}{\tan\sqrt\epsilon}\leq \frac{1}{\sqrt \epsilon}.\end{equation} Therefore for $E\in W_\epsilon^c$, $|\zeta_\epsilon(E)-1|\leq |(1+\epsilon)^{\kappa\gamma}-1|$.

  If $$c=\int_{\Gr_k(V)}f_1(E)dE\neq 0,$$ we find that $J_\epsilon\to c u(E_0)$ as $\epsilon\to 0$ by the dominated convergence theorem, while $|I_\epsilon|\to 0$ as $\epsilon\to 0$ by eq. \eqref{eq:I_epsilon}, and so for $\lambda_\epsilon=\epsilon^{-k(n-k)}$, $\lambda_\epsilon f_\epsilon \to c\delta_{E_0}$, concluding the proof.
 
 If $c=0$, we proceed as follows. We denote $T_0=T_{E_0}\Gr_k(V)$, and allow $\epsilon$ to assume arbitrary real values. Define $\tau:\Sigma^k_F\to T_0$ by $\tau(E)=\left.\frac{d}{d\epsilon}\right|_{\epsilon=0}g_\epsilon(E)\in T_0$. 
 
 \textit{Claim.} $\tau:\Sigma^k_F\to T_0$ is well-defined, and is in fact a diffeomorphism.
 
 \textit{Proof.} Let us first show that the claim holds on a
 neighborhood of $E_0$. Clearly $\tau(E_0)=0$. Let $(u_i)_{i=1}^k$ be an orthonormal basis of $E$. Using the Pl\"ucker embedding, we find that $$g_\epsilon(E)=\frac{\wedge_{i=1}^k (P_{E_0}u_i+\epsilon P_Fu_i)}{|\wedge_{i=1}^k (P_{E_0}u_i+\epsilon P_Fu_i)|}$$ is a smooth curve near $\epsilon=0$. Moreover, choosing the orthonormal basis $(u_i)$ smoothly in $E$ in a neighborhood of $E_0$, we see that $g_\epsilon(E)$ is smooth in $(\epsilon, E)$ in a neighborhhod of $(0, E_0)$, and so $\tau(E)=\left.\frac{d}{d\epsilon}\right|_{\epsilon=0}g_\epsilon(E)$ is well-defined and smooth near $E_0$.
 
 Let us now check that the derivative $D_{E_0}\tau$ is invertible. Choose any vectors $v_1,\dots, v_k\in F$, and consider the curve $E_t=\Span(e_1+tv_1,\dots,e_k+ tv_k)$ where $(e_i)$ is an orthonormal basis of $E_0$. Denote $w= \frac{d}{dt}|_{t=0}E_t$, which by \cite{kozlov13} is the general form of a tangent vector.
 
 We have $g_\epsilon(E_t)=\Span(e_1+t\epsilon v_1,\dots, e_k+t\epsilon v_k)$. As $|\wedge_{i=1}^k(e_i+t\epsilon v_i)|=1+O(\epsilon^2)$ for fixed $t$, we find
 \begin{equation}\label{eq:tau_formula} \tau(E_t)=\left.\frac{d}{d\epsilon}\right|_{\epsilon=0}\wedge_{i=1}^k (e_i+\epsilon t v_i)=t\sum_{i=1}^k e_1\wedge\dots\wedge e_{i-1}\wedge v_i\wedge e_{i+1}\wedge\dots\wedge e_k.\end{equation}
 Therefore  $$ D_{E_0}\tau (w)=\left.\frac{d}{dt}\right|_{t=0}\tau(E_t)=\sum_{i=1}^k e_1\wedge\dots\wedge e_{i-1}\wedge v_i\wedge e_{i+1}\wedge\dots\wedge  e_k.$$
 If $D_{E_0}\tau (w)=0$, we can wedge the sum separately with each $e_j$, $1\leq j\leq k$, and  conclude that $v_j=0$ for all $j$, so that $ w=0$. It follows that $D_{E_0}\tau$ has trivial kernel as required. 
 It now follows by the inverse mapping theorem that $\tau$ is a diffeomorphism of a neighborhood $U_0$ of $E_0$ onto a neighborhood $W_0$ of $0$.
 
 Observe that $g_s\circ g_t=g_{st}$ for all $s,t\in\R$, and so $\tau(g_s E)=s \tau(E)$ for all $s\in\R$ and $E\in\Sigma^k_F$. Note also that for any compact subset $K\subset \Sigma^k_F$, there is $\epsilon>0$ such that $g_\epsilon K\subset U_0$.
 As $\tau$ is smooth near $0$, it immediately follows that it is smooth on $\Sigma^k_F$. Similarly since $\tau$ is onto $W_0$, it is onto $T_0$. Finally if $E, E'\in \Sigma^k_F$ are distinct, choose $\epsilon>0$ such that $g_\epsilon E, g_\epsilon E'\in U_0$. Then $\tau(g_\epsilon E)\neq \tau(g_\epsilon E')$ and so $\tau(E)\neq \tau(E')$. Thus $\tau: \Sigma^k_F\to T_0$ is a diffeomorphism as claimed.
 
 We denote $b=\tau(E)$, $E\in \Sigma^k_F$.
 
 Given $E\in \Sigma^k_F$, one can fix an orthonormal basis $e_1,\dots, e_k$ of $E_0$ and orthonormal vectors $v_1,\dots, v_\kappa\in F$ such that $E=\Span(e_1+\cot \beta_1v_1,\dots, e_\kappa+\cot\beta_\kappa v_\kappa, e_{\kappa+1}, \dots, e_k)$.
 It then follows by eq. \eqref{eq:tau_formula} that \[b=\tau(E)=\sum_{i=1}^\kappa\cot\beta_i e_1\wedge\dots\wedge v_i\wedge\dots\wedge e_\kappa,\]
 and therefore
 \begin{equation}\label{eq:b_norm}|b|^2=\sum_{i=1}^\kappa \cot^2\beta_i.\end{equation}
  	
 	\textit{Claim.} $|\det D_b\tau^{-1}|$ is bounded from above on $T_0$. 
 	
 	\textit{Proof.} Indeed, differentiating the equation $\tau(g_s E)=s \tau(E)$ for $s>0$ fixed we find $D_E\tau=\frac{1}{s}\circ D_{g_sE}\tau\circ D_Eg_s$, where $D_Eg_s:T_E\Gr_k(V)\to  T_{g_sE}\Gr_k(V)$.
 	Taking determinants, we deduce that $\det D_E\tau=s^{-k(n-k)}\det ( D_{g_sE}\tau)\det (D_Eg_s)$.
 
 	Recall that $\Dens(T_E\Gr_k(V))=\Dens^{-n}(E)\otimes \Dens^k(V)$, so that \[\det (D_Eg_s)=\det(g_s)^k\Jac(g:E\to gE)^{-n}=s^{k(n-k)}\eta_{g_s}(E)^n,\]
and we conclude that \begin{equation}\label{eq:det_identity}\det D_E\tau= \eta_{g_s}(E)^n\det ( D_{g_sE}\tau).\end{equation}

Fix a compact neighborhood $K_0$ of $E_0$ inside $\Sigma^k_F$. Define $K_j=g_{2^j}K_0$, so that $K_0\subset K_1\subset \dots$ is an increasing sequence of compacts exhausting $\Sigma_F^k$. 

By eq. \eqref{eq:eta}, $\eta_{g_s}(E)\to 0$ as $s\to \infty$, uniformly on $K_1\setminus K_0$. As $K_{j+1}\setminus K_j =g_{2^j}(K_1\setminus K_0)$, it follows from eq. \eqref{eq:det_identity} that $|\det D_E\tau|\to\infty$ as $E\to \Xi^k_F$.
The claim readily follows.
 	 	
 	 	 It follows from eq. \eqref{eq:b_norm} that $|b|^2\leq \kappa \cot^2(\min\beta_i)$, and so by eq. \eqref{eq:f_1_first_bound} 
 	 	\begin{equation}\label{eq:f_1_bound}|f_1(E)|\leq c_0' e^{-c_1'/\min\beta_i}\leq c_0'e^{-c_1'\cot(\min\beta_i)}\leq c_0'e^{-c_1'|b|/\sqrt\kappa}.\end{equation}
 	 	
 	Write $\widetilde f_1(b)=f_1(\tau^{-1}b)|\det D_b\tau^{-1}|$, so that $\widetilde f_1(b)db=f_1(E)dE$. 
 	It follows from the boundedness of $|\det D_b\tau^{-1}|$ and eq. \eqref{eq:f_1_bound} that there are constants $\tilde c_0, \tilde c_1>0$ such that 
 	\begin{equation}\label{eq:tildef1_exp}|\widetilde f_1(b)| \leq \tilde c_0 e^{-\tilde c_1|b|}.\end{equation}
	 
 We see from eq. \eqref{eq:zeta} that $\zeta_{\epsilon}(E)$ is smooth in $(\epsilon,E)\in \R\times \Sigma^k_F$. Write its Taylor-Maclaurin series in $\epsilon^2$ as $$\zeta_{\epsilon}(E)\sim \sum_{j=0}^\infty \zeta_j(b)\epsilon^{2j}, \qquad \zeta_0(b)=1.$$
 By eqs. \eqref{eq:zeta} and \eqref{eq:b_norm}, there are constants $c_j$ such that  \begin{equation}\label{eq:zeta_bound}|\zeta_j(b)|\leq c_j|b|^{2j}.\end{equation} 
 Note that for $b\in \tau(W_\epsilon^c)$ it holds by eqs. \eqref{eq:cot_beta} and \eqref{eq:b_norm} that $|b|^2\leq \frac{\kappa}{\epsilon}$. Therefore $\epsilon|b|\leq \sqrt\kappa\epsilon^\frac 12$. 
 In particular for $b\in \tau(W_\epsilon^c)$, $\epsilon b\to 0$ as $\epsilon\to 0$ uniformly in  $b\in \tau(W_\epsilon^c)$, and so 
 $$\zeta_\epsilon(E)=\sum_{j=0}^N\zeta_j(b)\epsilon^{2j}+O(\epsilon^{2N+2}|b|^{2N+2}),\qquad   \epsilon\to0,\quad b\in\tau(W_\epsilon^c) .$$
 
 Since $\tau(g_\epsilon E)=\epsilon\tau(E)$, putting $\widetilde u(b)=u(\tau^{-1}b)$ we find $u(g_\epsilon E)=\widetilde u(\epsilon b)$. Denoting by $D^j\widetilde u(0)\in T_0^{\otimes j}$ the appropriate derivative of order $j$, write the Taylor-Maclaurin series in $\epsilon b$ as 
 $$\widetilde u(\epsilon b)\sim\sum_{j=0}^\infty \langle D^j\widetilde u(0), b^{\otimes j} \rangle\epsilon^{j},$$
 so that for $b\in \tau(W_\epsilon^c)$ we have, similarly to $\zeta_\epsilon$,  
$$\widetilde u(\epsilon b)=\sum_{j=0}^N\langle D^j\widetilde u(0), b^{\otimes j} \rangle\epsilon^{j}+O(\epsilon^{N+1}|b|^{N+1}),\qquad   \epsilon\to0,\quad b\in\tau(W_\epsilon^c) .$$
 
 We now make the change of variables $b=\tau (E)$ in eq. \eqref{eq:J_epsilon0}. Using the two approximations above for $\widetilde u(\epsilon b)$ and $\zeta_\epsilon(E)$, where the error terms are uniform in $b\in\tau(W_\epsilon^c)$, we dedude that for every $N\geq 0$,
 \begin{align}\label{eq:J_epsilon}\nonumber &J_\epsilon=\\&=\int_{\tau(W_\epsilon^c)} \widetilde f_1(b) (\sum_{j=0}^{N/2} \zeta_j(b)\epsilon^{2j} +O(\epsilon^{N+1}|b|^{N+1}))(\sum_{j=0}^N \langle D^j\widetilde u(0), b^{\otimes j} \rangle\epsilon^{j}+O(\epsilon^{N+1}|b|^{N+1})db\nonumber
 	\\&=\sum_{j=0}^N \langle \xi_j(\epsilon),\widetilde u\rangle\epsilon^j +O(\epsilon^{N+1}),\end{align}
 where each $\xi_j(\epsilon)$ is a distribution supported at $0$. Explicitly, 
 $$\langle \xi_j(\epsilon),\widetilde u\rangle=\langle D^j\widetilde u(0),\int_{\tau(W_\epsilon^c)} \widetilde f_1(b)b^{\otimes j}db\rangle +\sum_{i=1}^{j/2}\langle D^{j-2i}\widetilde u(0), \int_{\tau(W_\epsilon^c)}\widetilde f_1(b)\zeta_i(b) b^{\otimes (j-2i)} db\rangle .$$
  Introducing the distributions $\xi_j$ given by 
 $$\langle \xi_j,\widetilde u\rangle=\langle D^j\widetilde u(0),\int_{T_0} \widetilde f_1(b)b^{\otimes j}db\rangle +\sum_{i=1}^{j/2}\langle D^{j-2i}\widetilde u(0), \int_{T_0}\widetilde f_1(b)\zeta_i(b) b^{\otimes (j-2i)} db\rangle ,$$
 it follows from eqs. \eqref{eq:tildef1_exp} and \eqref{eq:zeta_bound} that for all $N$ there is $C_{N, j}$ such that
 \begin{equation}\label{eq:uniform_bound}\langle \xi_j(\epsilon)-\xi_j,\widetilde u\rangle|\leq C_{N,j}\|\widetilde u\|_{C^j(0)}\epsilon^N,\end{equation}
 where $\|\widetilde u\|_{C^j(0)}$ is the maximum of all derivatives at $0$ of order at most $j$.  In particular, the distributions $\xi_j(\epsilon)$ weakly converge to $\xi_j$ as $\epsilon\to 0$.
 
 Recall that by Borel's theorem, $D^j\widetilde u(0)$ can be arbitrary. If $\xi_j=0$ for all $j$, it must hold in particular that for all $j$, $\int_{T_0} \widetilde f_1(b)b^{\otimes j}db=0.$
 
 That is $\widetilde f_1(b)$, which is a smooth function on $T_0$, satisfies $\int_{T_0} \widetilde f_1(b)p(b)db=0$ for any polynomial $p$. Furthermore, by eq. \eqref{eq:tildef1_exp} it decays exponentially. Therefore, its Fourier transform $\mathcal F \widetilde f_1$ extends analytically to a neighborhood of the real subspace, while all its derivatives at the origin vanish. Consequently, $\widetilde f_1=0$. But $\widetilde f_1(b)=f_1(\tau^{-1}b)|\det D_b\tau^{-1}|$ is not identically zero, a contradiction.
 
 Thus there exists a least integer $m$ such that $\xi_m\neq 0$. It follows by eqs. \eqref{eq:J_epsilon} and \eqref{eq:uniform_bound} that
 
 $$ \epsilon^{-k(n-k)-m}\langle f_\epsilon, u\rangle \to \langle \xi_m, \widetilde u\rangle=\langle \tau^*\xi_m, u\rangle,\quad \epsilon\to 0^+.$$
 It remains to choose $\lambda_\epsilon=\epsilon^{-k(n-k)-m}$, so that $\lambda_\epsilon g_\epsilon\phi$ weakly converges to the generalized section $\phi_0$ which, under the Euclidean trivialization, corresponds to $\tau^*\xi_m$. 

\endproof

\begin{Theorem}\label{thm:unified}
	Let $W\subset \Gamma^{-\infty}(\Gr_k(V), M_\alpha)$ be a $\GL(V)$-invariant closed subspace such that its $\OO(n)$-types form a set which is not co-sparse. Then $W$ is Bernstein $\Xi$-quasianalytic, and $W\cap\Gamma^\infty(\Gr_k(V), M_\alpha)$ is exponentially $\Xi$-quasianalytic. 
\end{Theorem}
\proof
First we ought to show that a nonzero smooth section $\phi\in W$ cannot vanish at $\Xi^k$ exponentially.
Let $\phi$ vanish at $\Xi^k$ exponentially. Proposition \ref{prop:point_rescale} then guarantees the existence of a non-trivial $\phi_0\in \Gamma^{-\infty}(\Gr_{k}(V), M_\alpha)$ supported at a point, which still lies in $W$. This now contradicts Theorem \ref{thm:sparse}.

Now let $\phi$ be any generalized section in $W$, and assume $\Supp \phi\subset\Sigma^k$. Choose an approximate identity $\mu_\epsilon\in \mathcal M^\infty(\GL(V))$. Then $\mu_\epsilon\ast \phi$ is smooth, belongs to $W$, and $\Supp(\mu_\epsilon\ast \phi)\subset\Sigma^k$ for small enough $\epsilon$. It follows by the previous case that $\mu_\epsilon\ast \phi=0$, and taking $\epsilon\to 0$ we conclude that $\phi=0$.

\endproof

\begin{proof}[Proof of Theorem \ref{thm:alpha_cosine}, i).]
Follows at once from Theorem \ref{thm:unified} and Lemma \ref{lem:cosine_image} by taking $W=\mathrm{Image}(S_\alpha)$, and replacing $S_\alpha h$ with the corresponding section $\phi\in \Gamma^{-\infty}(\Gr_{n-k}(\R^n), M_\alpha)$ using the Euclidean structure and the orthogonal complement map $\Gr_k(\R^n)\to \Gr_{n-k}(\R^n)$. Note that the latter interchanges $\Sigma^k$ and $\Sigma^{n-k}$. 
\end{proof}

\begin{Remark}\label{rmk:gln_modules}
	Theorem \ref{thm:unified} can equally be applied to any of the proper submodules in the composition series of $\GL_n(\R)$ appearing in \cite{howe_lee_gln} in Theorems 3.4.2(ii) and 3.4.4(a)(ii). Those examples are particularly interesting, as it is unknown whether they lie in the kernel of an invariant differential operator.
\end{Remark}

\begin{proof}[Proof of Theorem \ref{thm:radon}.]
By taking orthogonal complements if necessary, we may assume $p<k$. Assume that either $h\in C^\infty(\Gr_p(\R^n))$ and $\mathcal R_{p,k}h$ vanishes at $\Xi^k$ exponentially, or $h\in C^{-\infty}(\Gr_p(\R^n))$ and $\Supp(\mathcal R_{p,k}h)\subset\Sigma^k$. Using Theorem \ref{thm:unified} and Lemma \ref{lem:radon_image}, we conclude that $\mathcal R_{p,k}h=0$ as in the proof of Theorem \ref{thm:alpha_cosine}. By Theorem \ref{thm:radon_injective}, $h=0$.
\end{proof}

\begin{Remark}
	An alternative way to conclude the proof exploits the fact that the range of the Radon transform coincides with the kernel of an invariant differential operator. Denoting this operator by $D_0$, as in Theorem \ref{thm:unified} we deduce the existence of $f_0=D_1\delta_{E_0}$ satisfying $D_0f_0=0$, for some $D_1\in U(\mathfrak{so}(n))$. Thus $D:=D_0\circ D_1\in U(\mathfrak{so}(n))$ satisfies $D\delta_{E_0}=0$. Proceeding as in the proof of Theorem \ref{thm:sparse}, we deduce the existence of a non-zero differential operator with polynomial coefficients which annihilates all polynomials.
\end{Remark}

As a corollary, we sharpen a theorem of Funk on sections of star bodies \cite{gardner}.
\begin{Corollary}\label{cor:funk_section}
	Assume $2\leq k\leq n-2$, and $S$ is a centrally-symmetric star body around $0\in\mathrm{int}(S)$. Let $U\subset\Gr_k(\R^n)$ be any open neighborhood of $\Xi^k$. Then $S$ is uniquely determined by its section function $A_S(E)=\vol_k(E\cap S)$, $E\in U$.
\end{Corollary}
\proof
Let $\rho_S:S^{n-1}\to \R$ be the radial function of $S$. Then $$A_S(E)=\frac{1}{k}\int_{E\cap S^{n-1}}\rho_S(\theta)^{k}d\theta=\frac1 k\mathcal R_{1, k}(\rho_S^{k}).$$ Theorem \ref{thm:radon} now concludes the proof.
\endproof

An immediate corollary of Theorem \ref{thm:cosine} is Theorem \ref{thm:klain_sharper}, which we now prove in a form allowing for generalized valuations.
\begin{Corollary}\label{cor:klain_generalized}
	Assume $2\leq k\leq n-2$. Let $\phi\in\Val_k^{-\infty,+}(\R^n)$ be a generalized, even, $k$-homogeneous valuation, whose Klain section is supported on $\Sigma^k$. Then $\phi=0$.
\end{Corollary}
\proof
Choose an approximate identity $\mu_\epsilon\in\mathcal M^\infty(\GL_n(\R))$ with compact support shrinking to $\{\mathrm{Id}\}$. For small $\epsilon$, $\Kl(\mu_\epsilon\ast \phi)=\mu_\epsilon\ast \Kl(\phi)$ is also supported in $\Sigma^k$, while by \cite[Theorem 1.3]{alesker_bernstein} we know that $\Kl(\mu_\epsilon\ast \phi)\in \mathrm{Image}(\Cos)$. By Theorem \ref{thm:cosine} it follows that $\Kl(\mu_\epsilon\ast \phi)=0$ as $\epsilon\to 0$, and by Klain's injectivity theorem we conclude that $\mu_\epsilon\ast \phi=0$ as $\epsilon\to 0$.
\endproof

The corresponding sharpening of Aleksandrov's projection theorem \cite{gardner} follows at once.
\begin{proof}[Proof of Corollary \ref{cor:alexandrov}]
For the valuation $\phi=V(K[k],\bullet[n-k])$ it holds that $\Kl(\phi)={n\choose k}f_K$. By Corollary \ref{cor:klain_generalized}, $f_K$ is uniquely determined by $f_K|_U$. By Aleksandrov's projection theorem, $K$ is determined by $f_K$.
\end{proof}
 
We now deduce the sharpening of Schneider's injectivity theorem from the corresponding result for Klain's injectivity theorem. We will make use of the Alesker product on smooth valuations \cite{alesker_product}, which commutes with restriction to subspaces, see e.g. \cite[Theorem 3.13]{alesker_survey}. Furthermore, the Alesker product extends to a product of continuous and smooth valuations, which takes values in continuous valuation \cite[Proposition 8.1.2]{alesker_fourier}.
\begin{Corollary}\label{cor:schneider_sharper}
		Assume $1\leq k\leq n-3$. Let $\phi\in\Val_k^{-}(\R^n)$ be a continuous, odd, $k$-homogeneous valuation, whose Schneider section $\Sc_\phi\in \Gamma(\Gr_{k+1}(V), \Val_k^{-}(E))$ is supported inside $\Sigma^{k+1}$. Then $\phi=0$.
\end{Corollary}
\proof
Take $\psi\in\Val_1^{-,\infty}(\R^n)$, and consider the Alesker product $\eta=\phi\cdot \psi\in\Val_{k+1}^{+}(\R^n)$. Since Alesker product commutes with restrictions to subspaces, $\Kl_\eta$ is supported inside $\Sigma^{k+1}$. By Corollary \ref{cor:klain_generalized}, $\eta=0$. It follows that $\phi\cdot \psi=0$ for all $\psi=\psi_1\cdots\psi_{n-k}$ with $\psi_j\in\Val_1^{\infty}(\R^n), \psi_1\in\Val_1^{-,\infty}(\R^n)$. Invoking Alesker's irreducibility, we conclude that $\phi\cdot \psi=0$ for all $\psi\in \Val_{n-k}^{-,\infty}(\R^n)$. By Alesker-Poincar\'e duality, $\phi=0$.
\endproof

It remains to prove the second part of Theorem \ref{thm:alpha_cosine}, which is a straightforward reduction to known facts.

\begin{proof}[Proof of Theorem \ref{thm:alpha_cosine}, ii).]
We may assume $k\leq \frac n 2$. It suffices to produce in each case a non-zero distribution in $\mathrm{Image}(S_\alpha)$ supported at $E_0$. Indeed, using the $\GL$-equivariant form of $S_\alpha$, convolving with an approximate identity on $\GL_n(\R)$ would then produce the desired smooth function of small support. 

When $S_\alpha$ is invertible (on smooth functions, or equivalently on distributions), the statement is trivial. This is the case if $\alpha\notin\mathbb Z$ by \cite[Theorem 4.15]{alesker_cosine}.

Similarly when $\alpha=-k$, by \cite[Theorem 1.6]{alesker_gourevitch_sahi} $S_\alpha$ is the Radon transform on $\Gr_k(\R^n)$ and is therefore invertible by \cite{gelfand_graev_rosu}. 

If $\alpha=-k-2m$ with $m\in\mathbb N$ then by Theorem \ref{thm:cosine_diff_operator} one can write $S_\alpha=DS_{-k}$ for some $\OO(n)$-equivariant differential operator $D$. Thus $D\delta_0$ lies in $\mathrm{\Image}(S_{\alpha})$ and is non-zero, or else $DZ_\lambda=0$ for all zonal harmonics $Z_\lambda$ at $E_0$ for all $\lambda\in\Lambda_k$, which by $\OO(n)$-equivariance would imply $D=0$, which is false. 

If $\alpha=-(k+1)-2m$ and $m\in\mathbb Z_+$, then again by Theorem \ref{thm:cosine_diff_operator} we can write $S_\alpha=D'S_{-(k+1)}$ for some $\OO(n)$-equivariant differential operator $D'$. 
If $n\neq 2k$ then by \cite[Corollary 1.6]{gourevitch_gln} we know that $S_{-(k+1)}$ is invertible, and we conclude as before.

Finally assume $n=2k$ and $\alpha=-(k+1)-2m$. Assume first that $m=0$. By \cite{howe_lee_gln}, $\Gamma(\Gr_{k}(\R^n), M_{-(k+1)})$ has composition series of length $2$, and we can find an irreducible submodule $W$ whose $\OO(n)$-types consist of those $\lambda\in \Lambda_k$ for which $\lambda_k>0$. Denote by $W^\perp\subset C^{-\infty}(\Gr_k(\R^n))$ the closed subspace generated by all remaining types, namely those with $\lambda_k=0$.
By Theorem \ref{thm:radon_range_pde} and Proposition \ref{prop:radon_types}, $\Ker\Omega_{k-1,k}=W^\perp$, while 
by \cite{gourevitch_gln}, $\mathrm{Image}(S_{-(k+1)})=W$. It follows that $f_0:=\Omega_{k-1,k}\delta_{E_0}\in\mathrm{Image}(S_{-(k+1)})$ is the sought-after distribution.

If $m\geq 1$ then as before we have $S_{-(k+1)-2m}=D'S_{-(k+1)}$. We know that $D'f_0=D'\Omega_{k-1,k}\delta_0\neq 0$, or else by the $\OO(n)$-equivariance of $D',\Omega_{k-1,k}$ it would follow that $D'\Omega_{k-1,k}=0$, which is impossible for invariant differential operators since the product of their principal symbols is nonzero.
Thus $D'f_0\in \mathrm{Image}(S_{-(k+1)-2m})$ is the desired distribution supported at $E_0$. This concludes the proof.
\end{proof}

We remark that one can easily avoid using \cite{gourevitch_gln}. In the case $n=2k$, $\alpha=-k-1$, the composition series have length $2$ by \cite{howe_lee_gln}, and so the two possibilities are these: either $S_{-(k+1)}$ is an isomorphism, or $\Image(S_{-(k+1)})=W$. Similarly for $n\neq 2k, \alpha=-k-1$, it suffices to note that the target space is irreducible \cite{howe_lee_gln} .

\section{Polynomial solutions of a PDE with polynomial coefficients}\label{sec:PDE}

Here we prove a general statement of independent interest about PDEs with polynomial coefficients, which roughly asserts that a partial differential equation with polynomial coefficients cannot have too many polynomial solutions. 
It is an essential ingredient in the proof of Theorem  \ref{thm:sparse}. Its proof was explained to us by Joseph Bernstein. 
We remark that in a previous version of this note, a somewhat more elementary though more involved proof was given for a weaker statement, yielding a correspondingly weaker form of Theorem \ref{thm:sparse}, asserting positive density for the support of the spectrum of a distribution supported at a point; this however suffices for the applications in this note. This is explained in Appendix \ref{app:PDE}.

Recall that $\mathcal P_m\subset \mathbb C[x_1,\dots, x_k]$ denotes the space of polynomials of degree at most $m$, and by Lemma \ref{lem:polynomials_dimension}, $\dim \mathcal P_m=\frac{1}{k!}m^k+O(m^{k-1})$ as $m\to\infty$.

\begin{Theorem}\label{thm:PDE_main}
	Let $D$ be a nonzero linear differential operator on $\mathbb C[x_1,\dots, x_k]$ with polynomial coefficients. Then $$ \dim(\Ker D\cap \mathcal P_m)=O(m^{k-1}), \quad m\to\infty.$$
\end{Theorem}
\proof
Identify the dual space $\mathcal P_m^*$ with the space $C^{-\infty,m}_{\{0\}}(\C^k)$ of distributions of order at most $m$ supported at $0$, with the natural pairing. Write also $\mathcal P_{-1}^*=\{0\}$. 
We may assume that the origin is placed in such a way that the principal symbol $\sigma_p(D)|_0\neq 0$. Let $N$ be the degree of $D$, and fix $M>0$ sufficiently large such that $D$ restricts to an operator $D_m:\mathcal P_m\to \mathcal P_{m+M}$. Its adjoint $D_m^*: \mathcal P_{m+M}^*\to \mathcal P_m^*$ is then given by a differential operator $D^*$ of order $N$ when restricted to $\mathcal P_{m-N}^*$. 

As $\sigma_p(D^*)|_0\neq 0$, we may write $D^*=L+D_0$, where $L$ is a differential operator of order $N$ with constant coefficients, and $\sigma_p(D_0)|_0=0$. 
In particular, $L$ and $D^*$ induce the same operator $T$ of order $N$ on the associated graded space, $T:\oplus_{i=0}^\infty \mathcal P_{i}^*/\mathcal P_{i-1}^*\to \oplus_{i=N}^\infty \mathcal P_{i}^*/\mathcal P_{i-1}^*.$ 

Recall that the action of $L$ on $C^{-\infty}_{\{0\}}(\C^k)$ is injective, e.g. because each $f\in C^{-\infty}_{\{0\}}(\C^k)$ can be uniquely represented as $f=L'\delta_0$ where $L'$ has constant coefficients \cite{hormander_vol1}. It follows that $T$ is injective, and therefore 
$D^*: C^{-\infty}_{\{0\}}(\C^k)\to C^{-\infty}_{\{0\}}(\C^k)$ is injective as well. Consequently,
\begin{align*} \dim \mathcal P_m-\dim \Ker D_m& =\dim(\Ker D_m)^\perp=\dim \mathrm{Image}(D_m^*)\\&\geq \dim(D^*(\mathcal P_{m-N}^*))=\dim \mathcal P_{m-N},\end{align*}
and therefore $\dim\Ker D_m\leq \dim\mathcal P_m-\dim\mathcal P_{m-N}=O(m^{k-1})$.

\endproof

\section{Examples and Discussion}\label{sec:discussion}

\subsection{Some examples of functions in the image of $\Cos$ with proper support}

It follows from \cite{alesker_bernstein} that any $\SO(n-1)$-invariant function or distribution on $\Gr_k(\R^n)$ belongs to the image of the cosine transform. In particular, one can find examples supported on the embedded grassmannians $\{E\subset \R^{n-1}\}$, resp. $\{E\supset \R^1\}$, or smooth examples with support in their respective arbitrarily small neighborhoods.

More examples can be obtained from the valuation theory of indefinite orthogonal groups, which was studied extensively in \cite{alesker_faifman_lorentz, bernig_faifman_opq}. In the space $\R^{p,q}$ with indefinite quadratic form $Q$, the space of $\OO(Q)$-invariant generalized valuations is $2$-dimensional. It is spanned by $\phi_k^0$, $\phi_k^1$, where $\Kl(\phi_k^i)$ is supported on the closure $S^i$ of the set of those subspaces $E$ where $Q|_E$ has signature $(p(E), q(E))$ with $q(E)\equiv i (\text{mod } 2)$. While those Klain sections are merely continuous, one can obtain smooth Klain sections with arbitrarily little change to the support by convolving with an approximate identity on $\GL_n(\R)$. By \cite[Theorem 1.3]{alesker_bernstein}, those sections lie in $\mathrm{Image}(\Cos)$.

Considering for example the $2$-homogeneous $\OO(2,2)$-invariant valuation $\phi_2^1$ on $\R^4$ with quadratic form $Q_\epsilon=x_1^2+x_2^2-\epsilon(x_3^2+x_4^2)$, and letting $\epsilon\to 0$, one arrives at a distribution in $\mathrm{Image}(\Cos)$ supported on $\{E: E\cap F\neq \{0\}\}$, where $F=\Span(e_3,e_4)$. Similarly, one can obtain elements in $\mathrm{Image}(\Cos)$ supported on $\{E\in \Gr_k(\R^n): E\cap F\neq \{0\}\}$ for any $2\leq j\leq n-2$ and $F\in\Gr_j(\R^n)$.

In light of these examples, we are led to consider the following more concrete version of Question \ref{qu:cosine}.
\begin{Question}
	Given a Schubert variety $X\subset \Gr_k(\R^n)$, is there $h\in C^{-\infty}(\Gr_k(\R^n))$ such that $\Cos(h)$ is supported on $X$? 
\end{Question}
Theorem \ref{thm:cosine} implies that the answer is negative when $X$ is a single point.

\subsection{Odd valuations}
For $F\in\Gr_{n-k}(\R^n)$, let $\Xi^{k+1}_F\subset \Gr_{k+1}(\R^n)$ be the collection of subsets that intersect $F$ non-generically, namely along a subspace of dimension $2$ or greater. It is easy to see that Theorem \ref{thm:klain_sharper} is equivalent to the following.
\begin{Theorem}\label{thm:klain_equivalent}
	Let $W\subset \Gr_{k+1}(\R^n)$ be a neighborhood of $\Xi^{k+1}_F$, and assume $2\leq k\leq n-2$. If $\phi\in\Val_k^+(\R^n)$ satisfies $\phi|_H=0$ for all $H\in W$ then $\phi=0$.
\end{Theorem}
\proof
Indeed, define $U\subset\Gr_k(\R^n)$ by $U=\cup_{H\in W}\Gr_k(H)$. Then $U$ is an open neighborhood of $\Xi^k_F$, and $\phi|_E=0$ for all $E\in U$. By Theorem $\ref{thm:klain_sharper}$, $\phi=0$.
\endproof
That Theorem \ref{thm:klain_sharper} is implied by Theorem \ref{thm:klain_equivalent} is similarly straightforward. Observe that unlike Theorem \ref{thm:klain_sharper}, this statement is not trivially false for odd valuations, and would constitute a strengthening of Corollary \ref{cor:schneider_sharper}.
\begin{Conjecture}
Theorem \ref{thm:klain_equivalent} holds also for odd valuations.
\end{Conjecture}

\subsection{On the support of distributions in the kernel of the cosine transform}
Complementing Question \ref{qu:cosine}, one could wonder about the possible supports of functions in the kernel of the cosine transform. In fact, there are no restrictions.

Indeed, $\Omega_{1,\kappa}$ annihilates all harmonics with $\lambda_2=0$ by Theorem \ref{thm:radon_range_pde} and Proposition \ref{prop:radon_types}, while the equation $S_{-1}=D_{-1,1}S_1$ of Theorem \ref{thm:cosine_diff_operator} implies that $D_{-1,1}$ annihilates all harmonics in $\mathrm{Image}(S_1)\cap \Ker(S_{-1})$, which by Proposition \ref{prop:cosine_types} includes the subset $\{\lambda\in\Lambda_\kappa: \lambda_2=2\}$. It follows that 
$D:=\Omega_{1,\kappa}\circ D_{-1,1}$ annihilates the types $\{\lambda\in\Lambda_\kappa:\lambda_2\leq 2\}$, which by Proposition \ref{prop:cosine_types} means that $\mathrm{Image}(S_1)\subset \Ker D$. Thus $f_E:=D\delta_{E}$ is supported at $\{E\}$, and is annihilated by $\Cos=S_1$. Taking Gelfand-Pettis integrals of $E\mapsto f_E$ with respect to a Borel measure on $\Gr_k(\R^n)$, one can obtain distributions in $\Ker(\Cos)$ of arbitrary support. Alternatively, passing to $\GL_n(\R)$-equivariant sections and convolving with an approximate identity on $\GL_n(\R)$, one can arrive at a smooth function in $\Ker\Cos$ with support in any given open set. 

\begin{Remark}
	It is a curious observation that $D$ annihilates all functions $f$ on $\Gr_k(\R^n)$ given by $f(E)=\vol_k(P_E(K))$ for some compact convex body $K$, as they lie in $\mathrm{Image}(\Cos)$. Note however that $\Ker D$ is strictly larger than $\mathrm{Image}(\Cos)$.
\end{Remark}
 
Observe that $f_E=D\delta_{E}$ is a distribution supported at a point with non co-finitely supported spectrum. Other such examples are $\Omega_{1,\kappa}\delta_{E}$ and $D_{-1,1}\delta_{E}$.

\subsection{Some further questions}
It would be interesting to consider distributions with small supports other than one point.
\begin{Question}\label{qu:small_support}
	Assume $f\in C^{-\infty}(\Gr_k(\R^n))$ is supported at a finite subset. Does it follow that 
	
	$$\liminf_{m\to\infty}\frac{|\{\lambda\in \Supp\widehat f:|\lambda|\leq2m\}|}{|\{\lambda\in \Lambda_\kappa:|\lambda|\leq2m\}|}>0?$$
	More generally, how small can $\Supp\widehat f$ get if $f$ is supported on a submanifold of given dimension?
\end{Question}

Finally, it would be interesting to study if one can strengthen Theorem \ref{thm:cosine} by replacing exponential quasianalyticity with the one of Denjoy-Carleman:

\begin{Question}
	Assume $2\leq k\leq n-2$, $h\in C^{\infty}(\Gr_k(\R^n))$, and $\Cos h$ vanishes on $\Xi^k$ with all derivatives of all orders. Does it follow that $\Cos h=0$? 
\end{Question}
One could then ask whether we can instead only assume the vanishing of derivatives of bounded order? At the other extreme, we have
\begin{Question}
	Does $\mathcal Ch|_{\Xi^k}=0$ with $h\in C^\infty(\Gr_k(\R^n))$ imply $\mathcal Ch=0$?
\end{Question} 
More generally, a likely difficult question with immediate applications to valuation theory and geometric tomography is describing the possible zero sets of functions in the image of the cosine transform. 

\appendix

\section{Positive density of polynomial solutions of a PDE}\label{app:PDE}
 	Here we prove a weaker version of Theorem \ref{thm:PDE_main} using elementary linear algebra.
 	
 	\begin{Theorem}\label{thm:PDE_weak}
 		Let $D$ be a nonzero linear differential operator on $\mathbb C[x_1,\dots, x_k]$ with polynomial coefficients. Then $$ \limsup_{m\to\infty}\frac {\dim(\Ker D\cap \mathcal P_m)}{\dim\mathcal P_m}<1,\quad m\to\infty.$$
 	\end{Theorem}
 	
 	It can be used instead of Theorem \ref{thm:PDE_main} in the proof of Theorem \ref{thm:sparse} to obtain the following weaker version of the latter:
 		\begin{Theorem}\label{thm:not_cosparse}
 		Assume $1\leq k\leq n-1$, and $f_0\in C^{-\infty}(\Gr_k(\R^n))$ is supported on a single $E_0\in \Gr_k(\R^n)$. Denote $\Lambda=\Supp\widehat {f_0}\subset \Lambda_{\kappa}$. Then $\Lambda$ is not sparse.
 	\end{Theorem}
 	Theorem \ref{thm:not_cosparse} can be used instead of Theorem \ref{thm:sparse} to establish a weaker form of Theorem \ref{thm:unified} with 'not co-sparse' replaced by 'sparse'. This then suffices to deduce Theorems \ref{thm:cosine}, \ref{thm:klain_sharper}, \ref{thm:alpha_cosine}, \ref{thm:radon}, since in all cases the $\OO(n)$-types appearing in the corresponding representations form in fact a sparse set by Lemmas \ref{lem:radon_image} and \ref{lem:cosine_image}.  	
 \newline\newline
 	We will need the following elementary determinant computation, which can also be deduced e.g. from \cite{yang_dong_determinants}.
 	
 	\begin{Lemma}\label{lem:factorial_matrix}
 		The $(n+1)\times (n+1)$ matrix $M=(\frac{1}{(k+i+j-2)!})_{i,j=1}^{n+1}$ 
 		$$ M=\begin{pmatrix}
 		 	\frac{1}{k!} & \frac{1}{(k+1)!} & \cdots & \frac{1}{(k+n)!}\\
 		 	\frac{1}{(k+1)!} & \frac{1}{(k+2)!} & \cdots & \frac{1}{(k+n+1)!} \\
 		 	\vdots  & \vdots  & \ddots & \vdots  \\
 		 	\frac{1}{(k+n)!} & \frac{1}{(k+n+1)!} & \cdots & \frac{1}{(k+2n)!}
 		 \end{pmatrix} $$
 		 is invertible for any integer $k\geq 0$. Moreover, 
 		
 		$$\det M=(-1)^{\frac{n(n+1)}{2}}\prod_{i=0}^n\frac{i!}{(k+n+i)!}.$$
 	\end{Lemma}
 \proof
  Multiplying the $i$-th row by $(k+i-1)!$ for all $i$, we find that $\det M=\frac{1}{\prod_{i=0}^{n}(k+i)!}\det M'$ where 
 \[M'=
\begin{pmatrix}
	1 & \frac{1}{k+1} & \cdots & \frac{1}{(k+1)\cdots(k+n)}\\
	1 & \frac{1}{k+2} & \cdots & \frac{1}{(k+1)\cdots(k+n+1)} \\
	\vdots  & \vdots  & \ddots & \vdots  \\
	1 & \frac{1}{k+n+1} & \cdots & \frac{1}{(k+n+1)\cdots(k+2n)}
\end{pmatrix}.
\]

We have the partial fraction decomposition $$\prod_{\nu=1}^{j}\frac{1}{x+\nu}=\sum_{\nu=1}^{j} \frac{c_{j,\nu}}{x+\nu},$$
with nonzero numerical coefficients $c_{j,\nu}$, and $c_{j,j}=\frac{(-1)^{j-1}}{(j-1)!}$. We will apply this formula to $x=k, k+1,\dots$. 

It follows that $\det M'=\prod_{j=1}^n c_{j,j}\det \widetilde M$ where

 \[\widetilde M=
\begin{pmatrix}
1 & \frac{1}{k+1} & \cdots & \frac{1}{k+n}\\
	1 & \frac{1}{k+2} & \cdots & \frac{1}{k+n+1} \\
	\vdots  & \vdots  & \ddots & \vdots  \\
	1 & \frac{1}{k+n+1} & \cdots & \frac{1}{k+2n}
\end{pmatrix}.
\]
Substracting the $(i+1)$-st row from the $i$-th row for $i=1,\dots, n$ and expanding along the first column, we find that $\det \widetilde M=(-1)^n\det \widehat M$, where  $\widehat M$ is the $n\times n$ matrix
 \[\widehat M=
\begin{pmatrix}
	\frac{1}{k+1}-\frac{1}{k+2} & \frac{1}{k+2}-\frac{1}{k+3} & \cdots & \frac{1}{k+n}-\frac{1}{k+n+1}\\
	 \frac{1}{k+2}-\frac{1}{k+3} &  \frac{1}{k+3}-\frac{1}{k+4} & \cdots &  \frac{1}{k+n+1}-\frac{1}{k+n+2} \\
	\vdots  & \vdots  & \ddots & \vdots  \\
	\frac{1}{k+n}-\frac{1}{k+n+1} & \frac{1}{k+n+1}-\frac{1}{k+n+2} & \cdots &\frac{1}{k+2n-1}-\frac{1}{k+2n}
\end{pmatrix}.
\]
Sequentially adding to the $j$-th column the sum of all subsequent columns for $j=1,\dots, n-1$, we get 

\begin{align*}\det \widehat M&=\det \begin{pmatrix}
	\frac{n}{(k+1)(k+n+1)} & \frac{n-1}{(k+2)(k+n+1)} & \cdots & \frac{1}{(k+n)(k+n+1)}\\
	\frac{n}{(k+2)(k+n+2)} & \frac{n-1}{(k+3)(k+n+2)} & \cdots & \frac{1}{(k+n+1)(k+n+2)} \\
	\vdots  & \vdots  & \ddots & \vdots  \\
	\frac{n}{(k+n)(k+2n)} & \frac{n-1}{(k+n+1)(k+2n)} & \cdots & \frac{1}{(k+2n-1)(k+2n)}
\end{pmatrix}
\\&=\frac{n!}{(k+n+1)(k+n+2)\cdots(k+2n)}\det \begin{pmatrix}
	\frac{1}{k+1} & \frac{1}{k+2} & \cdots & \frac{1}{k+n}\\
	\frac{1}{k+2} & \frac{1}{k+3} & \cdots & \frac{1}{k+n+1} \\
	\vdots  & \vdots  & \ddots & \vdots  \\
	\frac{1}{k+n} & \frac{1}{k+n+1} & \cdots & \frac{1}{k+2n-1}
\end{pmatrix}.
\end{align*}
The remaining determinant is a special case of the Cauchy determinant $\det(\frac{1}{x_i-y_j})$ with $x_i=k+i-1$ and $y_j=-j$, which then evaluates to  $$\frac{\prod_{i=2}^n\prod_{j=1}^{i-1}(i-j)^2}{\prod_{i,j=1}^n(k+i+j-1)}=\prod_{i=0}^{n-1}\frac{ (i!)^2(k+i)!}{(k+n+i)!}.$$
Collecting the various factors together yields the claimed value of $\det M$, which is clearly nonzero and so $M$ is invertible.
 \endproof
 We remark that for $k=1$, the inverse of $M$ was computed explicity in \cite{habermann_factorial}.
 \\\\
For $\alpha,\beta\in \mathbb Z^k$, write $\alpha\geq \beta$ if $\alpha_i\geq \beta_i$ for all $i$, $\|\alpha\|_\infty=\max_i |\alpha_i| $, $|\alpha|=\sum_i |\alpha_i|$.

\begin{proof}[Proof of Theorem \ref{thm:PDE_weak}.] 
Let $\alpha, \beta\in\mathbb Z_+^k$ be multi-indices, and write $A_N=\{\alpha\in \mathbb Z_+^k: \|\alpha\|_\infty\leq N-1\}$. One can write $D=\sum_{\alpha\in A_N}\sum_{\beta} c_{\alpha\beta} x^\beta\partial^{\alpha}$ for some $N$. We may further assume that all non-zero coefficients $c_{\alpha\beta}$ have $\alpha\geq \beta$, for we can replace $D$ by $D'=\partial^\alpha D$ with arbitrary $\alpha$ as large as necessary, and $\Ker D\subset \Ker D'$.

By Lemma \ref{lem:polynomials_dimension}, $\dim\mathcal P_m\sim \frac{m^k}{k!}$ as $m\to\infty$.
In particular   
$\lim_{m\to\infty}\frac{\dim \mathcal P_m}{\dim\mathcal P_{m+1}}=1$ and so $\lim_{m\to\infty}\frac{\dim \mathcal P_{m+2N}-\dim \mathcal P_m}{\dim\mathcal P_{m}}=0$. We may therefore restrict our attention to $m=2Nm'$ with $m'\to \infty$. Namely, we will prove 
$$\limsup_{m'\to \infty} \frac{\dim (\Ker D\cap \mathcal P_{2Nm'})}{\dim \mathcal P_{2Nm'}}<1.$$

Fix $m'$ and set $m=2Nm'$. Take $f=\sum _{|\gamma|\leq m} a_\gamma x^\gamma\in\mathcal P_m$. Then $$Df=\sum_{|\gamma|\leq m} a_\gamma\gamma! \sum_{\alpha,\beta}c_{\alpha\beta}\frac{1}{(\gamma-\alpha)!}x^{\gamma+\beta-\alpha}.$$
By construction $D(\mathcal P_m)\subset\mathcal P_m$, and so $Df=0$ if and only if $(a_\gamma\gamma!)_{|\gamma|\leq m}$ is in the kernel of 
$M=(\mu_{\sigma\gamma})_{|\sigma|,|\gamma|\leq m}$ which is given by
\begin{equation}\label{eq:matrix_coeff} \mu_{\sigma\gamma}=\sum_{\alpha-\beta=\gamma-\sigma}c_{\alpha\beta}\frac{1}{(\gamma-\alpha)!}.\end{equation}

For $j=(j_1,\dots, j_k)\in \mathbb N^k$ consider the blocks $$I_j=I_{j_1,\dots, j_k}=[2(j_1-1)N, 2j_1N-1]\times\dots\times [2(j_k-1)N, 2j_kN-1]\cap  \mathbb Z_+^k\subset \mathbb Z_+^k.$$ Denote $J(m')=\{j=(j_1,\dots, j_k)\in \mathbb N^k: j_1+\dots+j_k\leq m'\}$, and note that for $j\in J(m')$ and  $\gamma\in I_j$, one has $|\gamma|\leq m$.
Clearly $I_j\cap I_{j'}=\emptyset$ if $j\neq j'$.

Define also $I'_{j_1,\dots, j_k}\subset I_{j_1,\dots, j_k}$ by 
$$I'_{j_1,\dots, j_k}=[2j_1N-N, 2j_1N-1]\times\dots\times [2j_kN-N, 2j_kN-1].$$
It is easy to see that if $\gamma\in I'_j$, then $\mu_{\sigma\gamma}\neq 0$ only if $\sigma\in I_j$.
Indeed, $\mu_{\sigma\gamma}\neq 0$ implies $c_{\alpha\beta}\neq 0$ for some $\alpha,\beta\in A_N$ with $\alpha-\beta=\gamma-\sigma$. 
Thus $0\leq \gamma_i- \sigma_i\leq \alpha_i\leq N-1$ for all $i$, and $\gamma\in I'_j$ implies $\sigma\in I_j$.

It follows that 
\begin{equation*}\mathrm{rank} M\geq \sum_{j\in J(m')}\mathrm{rank} (\mu_{\sigma\gamma})_{\sigma\in I_j, \gamma\in I'_j}.\end{equation*}

Let us show that each matrix in the sum is nonzero, so that $\mathrm{rank} M\geq |J(m')|$. 
Choose a pair $\alpha', \beta'$ such that $c_{\alpha'\beta'}\neq0$, and denote $\delta'=\alpha'-\beta'$.
Fix $j\in J(m')$. By assumption the vector $(c_{\alpha, \alpha-\delta'})_{\alpha\in A_N}$ is non-zero. Examining \eqref{eq:matrix_coeff}, it suffices to show that the vectors $(\frac{1}{(\gamma-\alpha)!})_{\alpha\in A_N}\in\R^{A_N}$, as $\gamma$ ranges over $I'_{j}$, constitute a basis of $\R^{A_N}$. Identifying this space with $(\R^N)^{\otimes k}$, we have $(\frac{1}{(\gamma-\alpha)!})_{\alpha\in A_N}=(\frac{1}{(\gamma_1-\nu)!})_{\nu=0}^{N-1}\otimes\dots\otimes (\frac{1}{(\gamma_k-\nu)!})_{\nu=0}^{N-1}$. It therefore suffices to show that 
 the vectors $(\frac{1}{(z-\nu)!})_{\nu=0}^{N-1}\in \R^N$, as $z$ ranges over $[2j_iN-N, 2j_iN-1]\cap \mathbb Z$ for any given $1\leq i\leq k$, constitute a basis of $\R^N$. But this follows at once from Lemma \ref{lem:factorial_matrix}.
 
 It remains to note that by Lemma \ref{lem:polynomials_dimension}, $|J(m')|= \dim (x_1\cdots x_k\mathcal P_{m'-k})\sim \frac{1}{k!}(m')^k$ as $m'\to\infty$.
As $\dim (\Ker D\cap \mathcal P_{2Nm'})=\dim \mathcal P_{2Nm'}-\mathrm{rank}M$, we conclude that
$$\limsup_{m'\to \infty} \frac{\dim (\Ker D\cap \mathcal P_{2Nm'})}{\dim \mathcal P_{2Nm'}}\leq 1-\lim_{m'\to\infty}\frac{\frac1 {k!}(m')^k}{\frac{1}{k!}(2Nm')^k}=1-\frac{1}{2^kN^k}.$$

\end{proof}

\bibliographystyle{abbrv}
\bibliography{./ref_papers}

\begin{bibdiv}
\begin{biblist}

\bib{aleksandrov}{article}{
      author={Aleksandrov, Aleksandr~D.},
       title={On the theory of mixed volumes. {II}. {N}ew inequalities between
  mixed volumes and their application},
        date={1937},
     journal={Mat. Sbornik N.S.},
      volume={2},
       pages={1205\ndash 1238},
}

\bib{Alesker:Irreducibility}{article}{
      author={Alesker, S.},
       title={Description of translation invariant valuations on convex sets
  with solution of {P}. {M}c{M}ullen's conjecture},
        date={2001},
        ISSN={1016-443X},
     journal={Geom. Funct. Anal.},
      volume={11},
      number={2},
       pages={244\ndash 272},
         url={https://doi.org/10.1007/PL00001675},
      review={\MR{1837364}},
}

\bib{alesker_product}{article}{
      author={Alesker, S.},
       title={The multiplicative structure on continuous polynomial
  valuations},
        date={2004},
        ISSN={1016-443X},
     journal={Geom. Funct. Anal.},
      volume={14},
      number={1},
       pages={1\ndash 26},
         url={https://doi.org/10.1007/s00039-004-0450-2},
      review={\MR{2053598}},
}

\bib{alesker_survey}{article}{
      author={Alesker, Semyon},
       title={Theory of valuations on manifolds: a survey},
        date={2007},
        ISSN={1016-443X},
     journal={Geom. Funct. Anal.},
      volume={17},
      number={4},
       pages={1321\ndash 1341},
         url={https://doi.org/10.1007/s00039-007-0631-x},
      review={\MR{2373020}},
}

\bib{alesker_fourier}{article}{
      author={Alesker, Semyon},
       title={A {F}ourier-type transform on translation-invariant valuations on
  convex sets},
        date={2011},
        ISSN={0021-2172},
     journal={Israel J. Math.},
      volume={181},
       pages={189\ndash 294},
         url={https://doi.org/10.1007/s11856-011-0008-6},
      review={\MR{2773042}},
}

\bib{alesker_cosine}{incollection}{
      author={Alesker, Semyon},
       title={The {$\alpha$}-cosine transform and intertwining integrals on
  real {G}rassmannians},
        date={2012},
   booktitle={Geometric aspects of functional analysis},
      series={Lecture Notes in Math.},
      volume={2050},
   publisher={Springer, Heidelberg},
       pages={1\ndash 21},
         url={https://doi.org/10.1007/978-3-642-29849-3_1},
      review={\MR{2985123}},
}

\bib{alesker_bernstein}{article}{
      author={Alesker, Semyon},
      author={Bernstein, Joseph},
       title={Range characterization of the cosine transform on higher
  {G}rassmannians},
        date={2004},
        ISSN={0001-8708},
     journal={Adv. Math.},
      volume={184},
      number={2},
       pages={367\ndash 379},
         url={https://doi.org/10.1016/S0001-8708(03)00149-X},
      review={\MR{2054020}},
}

\bib{alesker_faifman_lorentz}{article}{
      author={Alesker, Semyon},
      author={Faifman, Dmitry},
       title={Convex valuations invariant under the {L}orentz group},
        date={2014},
        ISSN={0022-040X},
     journal={J. Differential Geom.},
      volume={98},
      number={2},
       pages={183\ndash 236},
         url={http://projecteuclid.org/euclid.jdg/1406552249},
      review={\MR{3238311}},
}

\bib{alesker_fu}{book}{
      author={Alesker, Semyon},
      author={Fu, Joseph H.~G.},
       title={Integral geometry and valuations},
      series={Advanced Courses in Mathematics. CRM Barcelona},
   publisher={Birkh\"{a}user/Springer, Basel},
        date={2014},
        ISBN={978-3-0348-0873-6; 978-3-0348-0874-3},
        note={Lectures from the Advanced Course on Integral Geometry and
  Valuation Theory held at the Centre de Recerca Matem\`atica (CRM), Barcelona,
  September 6--10, 2010, Edited by Eduardo Gallego and Gil Solanes},
      review={\MR{3380549}},
}

\bib{alesker_gourevitch_sahi}{article}{
      author={Alesker, Semyon},
      author={Gourevitch, Dmitry},
      author={Sahi, Siddhartha},
       title={On an analytic description of the {$\alpha$}-cosine transform on
  real {G}rassmannians},
        date={2016},
        ISSN={0219-1997},
     journal={Commun. Contemp. Math.},
      volume={18},
      number={2},
       pages={1550025, 43},
         url={https://doi.org/10.1142/S021919971550025X},
      review={\MR{3461427}},
}

\bib{amit_olevskii}{article}{
      author={Amit, Tomer},
      author={Olevskii, Alexander},
       title={On the annihilation of thin sets},
        date={2017},
      eprint={arXiv:1711.04131},
}

\bib{beckner}{article}{
      author={Beckner, William},
       title={Inequalities in {F}ourier analysis},
        date={1975},
        ISSN={0003-486X},
     journal={Ann. of Math. (2)},
      volume={102},
      number={1},
       pages={159\ndash 182},
         url={https://doi.org/10.2307/1970980},
      review={\MR{385456}},
}

\bib{benedicks_finite_measure}{article}{
      author={Benedicks, Michael},
       title={On {F}ourier transforms of functions supported on sets of finite
  {L}ebesgue measure},
        date={1985},
        ISSN={0022-247X},
     journal={J. Math. Anal. Appl.},
      volume={106},
      number={1},
       pages={180\ndash 183},
         url={https://doi.org/10.1016/0022-247X(85)90140-4},
      review={\MR{780328}},
}

\bib{bernig_faifman_opq}{article}{
      author={Bernig, Andreas},
      author={Faifman, Dmitry},
       title={Valuation theory of indefinite orthogonal groups},
        date={2017},
        ISSN={0022-1236},
     journal={J. Funct. Anal.},
      volume={273},
      number={6},
       pages={2167\ndash 2247},
         url={https://doi.org/10.1016/j.jfa.2017.06.005},
      review={\MR{3669033}},
}

\bib{bernstein_quasianalytic}{article}{
      author={Bernstein, Serge},
       title={Sur la d\'{e}finition et les propri\'{e}t\'{e}s des fonctions
  analytiques d'une variable r\'{e}elle},
        date={1914},
        ISSN={0025-5831},
     journal={Math. Ann.},
      volume={75},
      number={4},
       pages={449\ndash 468},
         url={https://doi.org/10.1007/BF01563654},
      review={\MR{1511806}},
}

\bib{beurling}{incollection}{
      author={Beurling, Arne},
       title={Sur les spectres des fonctions},
        date={1949},
   booktitle={Analyse {H}armonique},
      series={Colloques Internationaux du Centre National de la Recherche
  Scientifique, no. 15},
   publisher={Centre National de la Recherche Scientifique, Paris},
       pages={9\ndash 29},
      review={\MR{0033367}},
}

\bib{carleman_quasianalytic}{misc}{
      author={Carleman, T.},
       title={Les fonctions quasi analytiques. {Le{\c{c}}ons} profess{\'e}es au
  {Coll{\`e}ge} de {France}.},
    language={French},
         how={Collection de monographies sur la th{\'e}orie des fonctions.
  {Paris}: {Gauthier}-{Villars}. 115 p. (1926).},
        date={1926},
}

\bib{denjoy_quasianalytic}{article}{
      author={Denjoy, A.},
       title={Sur les fonctions quasi-analytiques de variable r{\'e}elle.},
    language={French},
        date={1921},
        ISSN={0001-4036},
     journal={C. R. Acad. Sci., Paris},
      volume={173},
       pages={1329\ndash 1331},
}

\bib{donoho_stark}{article}{
      author={Donoho, David~L.},
      author={Stark, Philip~B.},
       title={Uncertainty principles and signal recovery},
        date={1989},
        ISSN={0036-1399},
     journal={SIAM J. Appl. Math.},
      volume={49},
      number={3},
       pages={906\ndash 931},
         url={https://doi.org/10.1137/0149053},
      review={\MR{997928}},
}

\bib{folland_sitaram}{article}{
      author={Folland, Gerald~B.},
      author={Sitaram, Alladi},
       title={The uncertainty principle: a mathematical survey},
        date={1997},
        ISSN={1069-5869},
     journal={J. Fourier Anal. Appl.},
      volume={3},
      number={3},
       pages={207\ndash 238},
         url={https://doi.org/10.1007/BF02649110},
      review={\MR{1448337}},
}

\bib{gardner}{book}{
      author={Gardner, Richard~J.},
       title={Geometric tomography},
     edition={Second},
      series={Encyclopedia of Mathematics and its Applications},
   publisher={Cambridge University Press, New York},
        date={2006},
      volume={58},
        ISBN={0-521; 0-521-68493-5},
         url={https://doi.org/10.1017/CBO9781107341029},
      review={\MR{2251886}},
}

\bib{gelfand_graev_rosu}{article}{
      author={Gel'fand, I.~M.},
      author={Graev, M.~I.},
      author={Ro\c{s}u, R.},
       title={The problem of integral geometry and intertwining operators for a
  pair of real {G}rassmannian manifolds},
        date={1984},
        ISSN={0379-4024},
     journal={J. Operator Theory},
      volume={12},
      number={2},
       pages={359\ndash 383},
      review={\MR{757440}},
}

\bib{gonzalez_kakehi}{article}{
      author={Gonzalez, Fulton~B.},
      author={Kakehi, Tomoyuki},
       title={Pfaffian systems and {R}adon transforms on affine {G}rassmann
  manifolds},
        date={2003},
        ISSN={0025-5831},
     journal={Math. Ann.},
      volume={326},
      number={2},
       pages={237\ndash 273},
         url={https://doi.org/10.1007/s00208-002-0398-1},
      review={\MR{1990910}},
}

\bib{goodey_howard1}{article}{
      author={Goodey, Paul},
      author={Howard, Ralph},
       title={Processes of flats induced by higher-dimensional processes},
        date={1990},
        ISSN={0001-8708},
     journal={Adv. Math.},
      volume={80},
      number={1},
       pages={92\ndash 109},
         url={https://doi.org/10.1016/0001-8708(90)90016-G},
      review={\MR{1041885}},
}

\bib{goodey_schneider_weil_determination_projection}{article}{
      author={Goodey, Paul},
      author={Schneider, Rolf},
      author={Weil, Wolfgang},
       title={On the determination of convex bodies by projection functions},
        date={1997},
        ISSN={0024-6093},
     journal={Bull. London Math. Soc.},
      volume={29},
      number={1},
       pages={82\ndash 88},
         url={https://doi.org/10.1112/S0024609396001968},
      review={\MR{1416411}},
}

\bib{goodey_schneider_weil_directed}{article}{
      author={Goodey, Paul},
      author={Weil, Wolfgang},
       title={Directed projection functions of convex bodies},
        date={2006},
        ISSN={0026-9255},
     journal={Monatsh. Math.},
      volume={149},
      number={1},
       pages={43\ndash 64},
         url={https://doi.org/10.1007/s00605-005-0362-8},
      review={\MR{2260658}},
}

\bib{gourevitch_gln}{article}{
      author={Gourevitch, Dmitry},
       title={Composition series for degenerate principal series of {${\rm
  GL}(n)$}},
        date={2017},
        ISSN={0706-1994},
     journal={C. R. Math. Acad. Sci. Soc. R. Can.},
      volume={39},
      number={1},
       pages={1\ndash 12},
      review={\MR{3643436}},
}

\bib{grinberg_radon}{article}{
      author={Grinberg, Eric~L.},
       title={On images of {R}adon transforms},
        date={1985},
        ISSN={0012-7094},
     journal={Duke Math. J.},
      volume={52},
      number={4},
       pages={939\ndash 972},
         url={https://doi.org/10.1215/S0012-7094-85-05251-2},
      review={\MR{816395}},
}

\bib{grinberg_quinto}{article}{
      author={Grinberg, Eric~L.},
      author={Quinto, Eric~Todd},
       title={Analytic continuation of convex bodies and {F}unk's
  characterization of the sphere},
        date={2001},
        ISSN={0030-8730},
     journal={Pacific J. Math.},
      volume={201},
      number={2},
       pages={309\ndash 322},
         url={https://doi.org/10.2140/pjm.2001.201.309},
      review={\MR{1875896}},
}

\bib{groemer}{book}{
      author={Groemer, H.},
       title={Geometric applications of {F}ourier series and spherical
  harmonics},
      series={Encyclopedia of Mathematics and its Applications},
   publisher={Cambridge University Press, Cambridge},
        date={1996},
      volume={61},
        ISBN={0-521-47318-7},
         url={https://doi.org/10.1017/CBO9780511530005},
      review={\MR{1412143}},
}

\bib{habermann_factorial}{article}{
      author={Habermann, Karen},
       title={An explicit formula for the inverse of a factorial {H}ankel
  matrix},
        date={2021},
        ISSN={1034-4942},
     journal={Australas. J. Combin.},
      volume={79},
       pages={250\ndash 255},
      review={\MR{4224384}},
}

\bib{hadamard_question}{article}{
      author={Hadamard, J},
       title={Sur la g{\'e}n{\'e}ralisation de la notion de fonction
  analytique},
        date={1912},
     journal={Bull. Soc. Math. France},
      volume={40},
      number={suppl{\'e}ment sp{\'e}cial: vie de la soci{\'e}t{\'e}, s{\'e}ance
  du 28 f{\'e}vrier 1912},
       pages={28\ndash 29},
}

\bib{holmgren}{misc}{
      author={Holmgren, E.},
       title={{\"U}ber {Systeme} von linearen partiellen
  {Differentialgleichungen}.},
    language={German},
         how={Stockh. {\"O}fv. 58, 91-103 (1901).},
        date={1901},
}

\bib{hormander_vol1}{book}{
      author={H\"{o}rmander, Lars},
       title={The analysis of linear partial differential operators. {I}},
      series={Classics in Mathematics},
   publisher={Springer-Verlag, Berlin},
        date={2003},
        ISBN={3-540-00662-1},
         url={https://doi.org/10.1007/978-3-642-61497-2},
        note={Distribution theory and Fourier analysis, Reprint of the second
  (1990) edition [Springer, Berlin; MR1065993 (91m:35001a)]},
      review={\MR{1996773}},
}

\bib{howe_lee_gln}{article}{
      author={Howe, Roger},
      author={Lee, Soo~Teck},
       title={Degenerate principal series representations of {${\rm GL}_n(\bold
  C)$} and {${\rm GL}_n(\bold R)$}},
        date={1999},
        ISSN={0022-1236},
     journal={J. Funct. Anal.},
      volume={166},
      number={2},
       pages={244\ndash 309},
         url={https://doi.org/10.1006/jfan.1999.3427},
      review={\MR{1707754}},
}

\bib{james_constantine}{article}{
      author={James, Alan~T.},
      author={Constantine, A.~G.},
       title={Generalized {J}acobi polynomials as spherical functions of the
  {G}rassmann manifold},
        date={1974},
        ISSN={0024-6115},
     journal={Proc. London Math. Soc. (3)},
      volume={29},
       pages={174\ndash 192},
         url={https://doi.org/10.1112/plms/s3-29.1.174},
      review={\MR{374523}},
}

\bib{kakehi99}{article}{
      author={Kakehi, Tomoyuki},
       title={Integral geometry on {G}rassmann manifolds and calculus of
  invariant differential operators},
        date={1999},
        ISSN={0022-1236},
     journal={J. Funct. Anal.},
      volume={168},
      number={1},
       pages={1\ndash 45},
         url={https://doi.org/10.1006/jfan.1999.3459},
      review={\MR{1717855}},
}

\bib{kazdan}{incollection}{
      author={Kazdan, Jerry},
       title={Quasi-analytic functions},
        date={1968},
   booktitle={Entire functions and related parts of analysis},
      editor={Korevaar, Jacob},
      editor={Chern, S.~S.},
      editor={Ehrenpreis, Leon},
      editor={Fuchs, W. H.~J.},
      editor={Rubel, L.~A.},
   publisher={American Mathematical Society, Providence, R.I.},
       pages={244\ndash 252},
      review={\MR{0232935}},
}

\bib{klain_even}{article}{
      author={Klain, Daniel~A.},
       title={Even valuations on convex bodies},
        date={2000},
        ISSN={0002-9947},
     journal={Trans. Amer. Math. Soc.},
      volume={352},
      number={1},
       pages={71\ndash 93},
         url={https://doi.org/10.1090/S0002-9947-99-02240-0},
      review={\MR{1487620}},
}

\bib{kozlov_eng12}{article}{
      author={Kozlov, S.~E.},
       title={Geometry of real {Grassmannian} manifolds. {I}, {II}},
    language={English},
        date={1997},
        ISSN={1072-3374},
     journal={J. Math. Sci., New York},
      volume={100},
      number={3},
       pages={2239\ndash 2253},
}

\bib{kozlov13}{article}{
      author={Kozlov, S.~E.},
       title={Geometry of real {G}rassmannian manifolds. {I}, {II}, {III}},
        date={1997},
        ISSN={0373-2703},
     journal={Zap. Nauchn. Sem. S.-Peterburg. Otdel. Mat. Inst. Steklov.
  (POMI)},
      volume={246},
      number={Geom. i Topol. 2},
       pages={84\ndash 107, 108\ndash 129, 197\ndash 198},
         url={https://doi.org/10.1007/s10958-000-0008-2},
      review={\MR{1631784}},
}

\bib{kozlov_english_3}{article}{
      author={Kozlov, S.~E.},
       title={Geometry of real {Grassmannian} manifolds. {III}},
    language={English},
        date={1997},
        ISSN={1072-3374},
     journal={J. Math. Sci., New York},
      volume={100},
      number={3},
       pages={2254\ndash 2268},
}

\bib{matheron_book}{book}{
      author={Matheron, G.},
       title={Random sets and integral geometry},
      series={Wiley Series in Probability and Mathematical Statistics},
   publisher={John Wiley \& Sons, New York-London-Sydney},
        date={1975},
        note={With a foreword by Geoffrey S. Watson},
      review={\MR{0385969}},
}

\bib{mcmullen_decomposition}{article}{
      author={McMullen, P.},
       title={Valuations and {E}uler-type relations on certain classes of
  convex polytopes},
        date={1977},
        ISSN={0024-6115},
     journal={Proc. London Math. Soc. (3)},
      volume={35},
      number={1},
       pages={113\ndash 135},
         url={https://doi.org/10.1112/plms/s3-35.1.113},
      review={\MR{448239}},
}

\bib{meshulam06}{article}{
      author={Meshulam, Roy},
       title={An uncertainty inequality for finite abelian groups},
        date={2006},
        ISSN={0195-6698},
     journal={European J. Combin.},
      volume={27},
      number={1},
       pages={63\ndash 67},
         url={https://doi.org/10.1016/j.ejc.2004.07.009},
      review={\MR{2186416}},
}

\bib{narayanan_sitaram_lacunary}{article}{
      author={Narayanan, E.~K.},
      author={Sitaram, A.},
       title={Lacunary {F}ourier series and a qualitative uncertainty principle
  for compact {L}ie groups},
        date={2011},
        ISSN={0253-4142},
     journal={Proc. Indian Acad. Sci. Math. Sci.},
      volume={121},
      number={1},
       pages={77\ndash 81},
         url={https://doi.org/10.1007/s12044-011-0006-y},
      review={\MR{2814316}},
}

\bib{mandelbrojt_nazarov}{article}{
      author={Nazarov, Fedor},
       title={Local estimates for exponential polynomials and their
  applications to inequalities of the uncertainty principle type},
        date={1993},
        ISSN={0234-0852},
     journal={Algebra i Analiz},
      volume={5},
      number={4},
       pages={3\ndash 66},
      review={\MR{1246419}},
}

\bib{nazarov_olevskii}{incollection}{
      author={Nazarov, Fedor},
      author={Olevskii, Alexander},
       title={A function with support of finite measure and ``small''
  spectrum},
        date={2018},
   booktitle={50 years with {H}ardy spaces},
      series={Oper. Theory Adv. Appl.},
      volume={261},
   publisher={Birkh\"{a}user/Springer, Cham},
       pages={389\ndash 393},
      review={\MR{3792106}},
}

\bib{nazarov_ryabogin_zvavitch}{article}{
      author={Nazarov, Fedor},
      author={Ryabogin, Dmitry},
      author={Zvavitch, Artem},
       title={Non-uniqueness of convex bodies with prescribed volumes of
  sections and projections},
        date={2013},
        ISSN={0025-5793},
     journal={Mathematika},
      volume={59},
      number={1},
       pages={213\ndash 221},
         url={https://doi.org/10.1112/S0025579312001027},
      review={\MR{3028180}},
}

\bib{olafsson_cosine}{incollection}{
      author={\'{O}lafsson, G.},
      author={Pasquale, A.},
      author={Rubin, B.},
       title={Analytic and group-theoretic aspects of the cosine transform},
        date={2013},
   booktitle={Geometric analysis and integral geometry},
      series={Contemp. Math.},
      volume={598},
   publisher={Amer. Math. Soc., Providence, RI},
       pages={167\ndash 188},
         url={https://doi.org/10.1090/conm/598/12009},
      review={\MR{3156446}},
}

\bib{rubin_cosine}{article}{
      author={Rubin, B.},
       title={Funk, cosine, and sine transforms on {S}tiefel and {G}rassmann
  manifolds},
        date={2013},
        ISSN={1050-6926},
     journal={J. Geom. Anal.},
      volume={23},
      number={3},
       pages={1441\ndash 1497},
         url={https://doi.org/10.1007/s12220-012-9294-4},
      review={\MR{3078361}},
}

\bib{schneider_polytopes_projections}{article}{
      author={Schneider, Rolf},
       title={On the projections of a convex polytope},
        date={1970},
        ISSN={0030-8730},
     journal={Pacific J. Math.},
      volume={32},
       pages={799\ndash 803},
         url={http://projecteuclid.org/euclid.pjm/1117559014},
      review={\MR{267461}},
}

\bib{schenider_cosine}{article}{
      author={Schneider, Rolf},
       title={\"{U}ber eine {I}ntegralgleichung in der {T}heorie der konvexen
  {K}\"{o}rper},
        date={1970},
        ISSN={0025-584X},
     journal={Math. Nachr.},
      volume={44},
       pages={55\ndash 75},
         url={https://doi.org/10.1002/mana.19700440105},
      review={\MR{275286}},
}

\bib{schneider_simple}{article}{
      author={Schneider, Rolf},
       title={Simple valuations on convex bodies},
        date={1996},
        ISSN={0025-5793},
     journal={Mathematika},
      volume={43},
      number={1},
       pages={32\ndash 39},
         url={https://doi.org/10.1112/S0025579300011578},
      review={\MR{1401706}},
}

\bib{schneider_brightness}{article}{
      author={Schneider, Rolf},
       title={Verification of polytopes by brightness functions},
        date={2009},
        ISSN={0002-9939},
     journal={Proc. Amer. Math. Soc.},
      volume={137},
      number={11},
       pages={3899\ndash 3903},
         url={https://doi.org/10.1090/S0002-9939-09-10041-2},
      review={\MR{2529898}},
}

\bib{schneider_book}{book}{
      author={Schneider, Rolf},
       title={Convex bodies: the {B}runn-{M}inkowski theory},
     edition={expanded},
      series={Encyclopedia of Mathematics and its Applications},
   publisher={Cambridge University Press, Cambridge},
        date={2014},
      volume={151},
        ISBN={978-1-107-60101-7},
      review={\MR{3155183}},
}

\bib{schneider_weil_polytopes}{article}{
      author={Schneider, Rolf},
      author={Weil, Wolfgang},
       title={\"{U}ber die {B}estimmung eines konvexen {K}\"{o}rpers durch die
  {I}nhalte seiner {P}rojektionen},
        date={1970},
        ISSN={0025-5874},
     journal={Math. Z.},
      volume={116},
       pages={338\ndash 348},
         url={https://doi.org/10.1007/BF01111841},
      review={\MR{283692}},
}

\bib{sugiura}{article}{
      author={Sugiura, Mitsuo},
       title={Representations of compact groups realized by spherical functions
  on symmetric spaces},
        date={1962},
        ISSN={0021-4280},
     journal={Proc. Japan Acad.},
      volume={38},
       pages={111\ndash 113},
         url={http://projecteuclid.org/euclid.pja/1195523466},
      review={\MR{142689}},
}

\bib{szekeres_partitions}{article}{
      author={Szekeres, G.},
       title={An asymptotic formula in the theory of partitions},
        date={1951},
        ISSN={0033-5606},
     journal={Quart. J. Math. Oxford Ser. (2)},
      volume={2},
       pages={85\ndash 108},
         url={https://doi.org/10.1093/qmath/2.1.85},
      review={\MR{43129}},
}

\bib{tao05}{article}{
      author={Tao, Terence},
       title={An uncertainty principle for cyclic groups of prime order},
        date={2005},
        ISSN={1073-2780},
     journal={Math. Res. Lett.},
      volume={12},
      number={1},
       pages={121\ndash 127},
         url={https://doi.org/10.4310/MRL.2005.v12.n1.a11},
      review={\MR{2122735}},
}

\bib{zhelobenko}{book}{
      author={\v{Z}elobenko, D.~P.},
       title={Compact {L}ie groups and their representations},
      series={Translations of Mathematical Monographs, Vol. 40},
   publisher={American Mathematical Society, Providence, R.I.},
        date={1973},
        note={Translated from the Russian by Israel Program for Scientific
  Translations},
      review={\MR{0473098}},
}

\bib{wallach}{book}{
      author={Wallach, Nolan~R.},
       title={Real reductive groups. {II}},
      series={Pure and Applied Mathematics},
   publisher={Academic Press, Inc., Boston, MA},
        date={1992},
      volume={132},
        ISBN={0-12-732961-7},
      review={\MR{1170566}},
}

\bib{wigderson2}{article}{
      author={Wigderson, Avi},
      author={Wigderson, Yuval},
       title={The uncertainty principle: variations on a theme},
        date={2021},
        ISSN={0273-0979},
     journal={Bull. Amer. Math. Soc. (N.S.)},
      volume={58},
      number={2},
       pages={225\ndash 261},
         url={https://doi.org/10.1090/bull/1715},
      review={\MR{4229152}},
}

\bib{yang_dong_determinants}{article}{
      author={Yang, Sheng-Liang},
      author={Dong, Yan-Ni},
       title={Hankel determinants of the generalized factorials},
        date={2018},
        ISSN={0019-5588},
     journal={Indian J. Pure Appl. Math.},
      volume={49},
      number={2},
       pages={217\ndash 225},
         url={https://doi.org/10.1007/s13226-018-0264-9},
      review={\MR{3802784}},
}

\bib{zhang_cosine}{article}{
      author={Zhang, Genkai},
       title={Radon, cosine and sine transforms on {G}rassmannian manifolds},
        date={2009},
        ISSN={1073-7928},
     journal={Int. Math. Res. Not. IMRN},
      number={10},
       pages={1743\ndash 1772},
         url={https://doi.org/10.1093/imrn/rnn170},
      review={\MR{2505576}},
}

\bib{zygmund_trigonometric}{book}{
      author={Zygmund, A.},
       title={Trigonometric series. {V}ol. {I}, {II}},
     edition={Third},
      series={Cambridge Mathematical Library},
   publisher={Cambridge University Press, Cambridge},
        date={2002},
        ISBN={0-521-89053-5},
        note={With a foreword by Robert A. Fefferman},
      review={\MR{1963498}},
}

\end{biblist}
\end{bibdiv}

 \end{document}